\documentclass[centertags,leqno,11pt]{article}
\RequirePackage[mathscr]{eucal}
\RequirePackage{xspace,amssymb,mathrsfs,url,upref,verbatim}
\let\mathcal\mathscr
\usepackage{amsmath}
\usepackage{amssymb}
\usepackage{latexsym}
\usepackage{amsxtra}
\usepackage{amscd}
\usepackage{theorem}
\usepackage{fullpage}
\usepackage{titlesec}
\titleformat{\subsection}[runin]{\normalfont\bfseries}{\thesubsection.}{3pt}{}

\usepackage{graphics}
\usepackage{epic}

\renewcommand{\bar}{\overline}

\newcommand{\proof}{\noindent {\bf Proof:\ }}
\newcommand{\Endproof}{\hspace*{\fill} $\Box$ \vspace{1ex} \noindent }

\makeatletter
\renewcommand{\subsection}{\@startsection{subsection}{2}%
{\z@}{-3.25ex plus -1ex minus-.2ex}{-1em}{\bf}} \makeatother

\newcommand{\PP}{\mathbb{P}}
\newcommand{\ZZ}{\mathbb{Z}}

\newcommand{\QQ}{\mathbb{Q}}
\newcommand{\NN}{\mathbb{N}}
\newcommand{\FF}{\mathbb{F}}
\renewcommand{\AA}{\mathbb{A}}
\newcommand{\GG}{\mathbb{G}}

\newcommand{\tr}{{\rm trace}}

\newcommand{\bQl}{{\overline{\QQ}_\ell}}

\newcommand{\SL}{{\rm SL}}

\newcommand{\Gal}{{\rm Gal}}

\newcommand{\Spec}{{\rm Spec\,}}

\newcommand{\Hom}{{\rm Hom}}

\renewcommand{\middle}{{\rm mid}}

\newcommand{\MC}{{\rm MC}}

\newcommand{\To}{\;\longrightarrow\;}

\newcommand{\Mapsto}{\;\longmapsto\;}

\newcommand{\im}{{\rm Im}}

\newcommand{\pr}{{{\rm pr}}}

\newcommand{\rk}{{\rm rk}}

\newcommand{\Frob}{{\rm Frob}}

\newcommand{\Perv}{{\rm Perv}}

\renewcommand{\L}{{\mathcal L}}

\newcommand{\1}{{\bf 1}}

\newcommand{\cF}{{\mathcal F}}
\newcommand{\cG}{{\mathcal G}}
\newcommand{\cH}{{\mathcal H}}
\newcommand{\cJ}{{\mathcal J}}
\newcommand{\cK}{{\mathcal K}}
\newcommand{\cL}{{\mathcal L}}

\newcommand{\vD}{{\bf D}}
\newcommand{\rD}{{\rm D}}
\newcommand{\cP}{{\mathcal P}}
\newcommand{\Rhom}{{\rm Rhom}}

\newcommand{\Ind}{{\rm Ind}}
\newcommand{\ob}{{\rm ob}\,}

\newcommand{\Gr}{{\rm Gr}}

\theoremstyle{change}

{\theorembodyfont{\slshape}

\newtheorem{thm}{Theorem.}[section]
\newtheorem{cor}[thm]{Corollary.}
\newtheorem{lem}[thm]{Lemma.}
\newtheorem{prop}[thm]{Proposition.}

\newtheorem{defn}[thm]{Definition.} 
}

{\theorembodyfont{\rmfamily}

\newtheorem{rem}[thm]{Remark.}

}

\numberwithin{equation}{section}
\numberwithin{thm}{subsection}
\theoremstyle{plain}

\title{On the $\ell$-adic Fourier transform  and the determinant of the middle convolution}
\author{
   Michael Dettweiler
  }

\begin{document}
\maketitle 
\begin{abstract} 
We study the relation of the middle convolution to the $\ell$-adic Fourier transformation in the \'etale context. Using Katz' work and 
Laumon's theory of local Fourier transformations we obtain 
 a detailed description of the local monodromy and the determinant of Katz' middle convolution functor $\MC_\chi$ in the tame case.  
 The theory of local $\epsilon$-constants then implies that the property of an \'etale sheaf of having an  at most quadratic
 determinant  up to Tate twist is often preserved under $\MC_\chi$ if $\chi$ is quadratic. 
\end{abstract}

\tableofcontents 

%------------------------------------------------------

%\addcontentsline{toc}{section}{Introduction}
%\setcounter{page}{1}
\section*{Introduction}\label{Introduction}

Consider the addition map 
$ \pi: \AA^n_k\times \AA^n_k \to \AA^n_k$ for $k$  a finite  field. 
If $K$ and $L$ are objects  in the derived category  $\rD^b_c(\AA^1_k,\overline{\QQ}_\ell)$
  then one may consider two kinds of convolutions, exchanged by Verdier duality (cf.~\cite{Katz96}):
  $$ K\ast_* L=R\pi_*(K \boxtimes L) \quad \textrm{and}\quad K\ast_! L:=R\pi_!(K \boxtimes L).$$
   It is convenient  to restrict the above construction to the abelian category of perverse sheaves 
$\Perv(\AA^n,\overline{\QQ}_\ell) \subseteq \rD^b_c(\AA^n,\overline{\QQ}_\ell)$ (see~\cite{BBD}). Under some restrictions (e.g., if $n=1$ and if 
$K$ is geometrically irreducible and not geometrically translation invariant, \cite{Katz96}, Lem.~2.6.9)
the above defined convolutions are again perverse and  one can define the 
{\it middle convolution} of $K$ and $L\in \Perv(\AA^n,\overline{\QQ}_\ell)$ 
as 
\begin{equation}\label{eqimmid}
 K*_{{\rm mid}}L=\im\left(K\ast_! L \to K\ast_* L\right),\nonumber \end{equation}
cf.~\cite{Katz96}, Chap.~2.6. 

% We want to remark that although in many cases 
% the middle convolution can be expressed concretely in terms of 
% sheaf cohomology, avoiding the language of perverse sheaves, the basic properties  of the middle convolution, 
% like associativity,  can only be understood 
% in the larger framework of $\rD^b_c(\AA^1,\bQl)$ and perverse sheaves. \\
 
One reason why one is interested in the middle 
 convolution is that 
 $K\ast_{\rm mid} L$  is 
often irreducible, while  the convolutions
$K\ast_* L$ and $K\ast_! L$
 are usually mixed and hence not irreducible. 
 A striking application of  the concept of middle 
convolution is Katz' existence algorithm for irreducible 
rigid local systems (\cite{Katz96}, Chap.~6).\\

The aim of this article is the determination of the behaviour  of the Frobenius determinants under the middle convolution with Kummer sheaves.\\

 Our main results are: 
 \begin{enumerate}
\item Using Laumon's theory of local Fourier transformation \cite{Laumon} and the principle of stationary phase
(\cite{Laumon}, \cite{KatzTraveauxdeLaumon}) we derive in Thm.~\ref{thmThom} an explicit description of the local monodromy (the structure of 
Frobenius elements on the vanishing cycle spaces at the singularities) 
$$\MC_\chi(K):=K*_\middle (j_*\cL_\chi[1]),$$ for $K$ an irreducible nontrivial 
 tame middle 
extension sheaf and $\cL_\chi$ a Kummer sheaf. 
\item Building on Laumon's product formula expressing  the epsilon constant in terms of local epsilon constants \cite{Laumon}, we
obtain a formula for the determinant of $\MC_\chi(K)$ in the tame case (Cor.~\ref{cormcchidet}).
\item From Thm.~\ref{thmThom} and Cor.~\ref{cormcchidet} we conclude in Thm.~\ref{thml1} that, under certain natural restrictions, 
the property for a tame middle extension sheaf of having an 
at most quadratic determinant up to Tate-twist is preserved under middle convolution with quadratic Kummer sheaves.
 \end{enumerate}

   In a companion  article to this work, written jointly with Stefan Reiter \cite{DR23}, we use these methods in order to   
 prove the following: {\it Let $\FF_q$ be the finite field 
of order $q=\ell^k,$ where $\ell$ is an odd prime number and $k\in \NN .$ 
Then the special linear group $\SL_{n}(\FF_q)$  
occurs regularly as Galois
group over $\QQ(t)$ if $n>8\varphi(q-1)+11$ and if $q$ is odd.} 

%Another application of our methods is that they 
%allow to accompany the above mentioned algorithm of Katz for quasiunipotent rigid local systems 
%with an algorithm which gives the Frobenius traces (at smooth points and at the nearby cycle spaces at the singularities) 
%in each step. This enables the computation of the unramified  local $L$-functions associated to the Galois representations associated 
%to rigid local systems (cf.~\cite{DR07}, \cite{DettweilerSabbah}, \cite{PatrikisTaylor}).

%{\bf Acknowledgments:} This article presents a revised part of my 
%Habilitation Thesis (Heidelberg, 2005). 
%I  heartily thank  B.H. Matzat for 
%his support during the last years
% and  for many valuable discussions 
%on the subject of this work. 
%I thank M. Berkenbosch, D. Haran
%J. Hartmann, M. Jarden, U. K\"uhn,  S. Reiter, A. R\"oscheisen and S. Wewers 
%for valuable comments and discussions. 
%Part of this work was written 
%during my stays at the School of Mathematics of the Tel Aviv University
%(spring 2004), the Institute for Advanced Study (IAS)
%in Princetion (spring 2007), and at the Insitut des Hautes \'Etudes
%Scientifiques (IHES), Bures sur Yvette, Paris (Fall 2008). I thank these institutes 
%for creating a very friendly  and inspiring atmosphere. 

\section{General notation and conventions}\label{secvondefetale} 
 \subsection{General notation.} If $F$ is any field, then $\bar{F}$ denotes an algebraic closure of $F.$   Let in the following $k$ be  an either  finite or  algebraically 
 closed  field  of characteristic $p\geq 0$  and 
 let $\ell$ be a prime $\ell\neq p.$  
Let us  recall the setup used in Laumon's work on the $\ell$-adic Fourier transform,
see  \cite{Laumon} and the references therein:\\

If $X$ is a variety over $k$ (meaning that $X$ is separated of finite type over $k$), then $|X|$ denotes the set of closed points of $X.$
For $x\in |X|,$ the residual field is denoted $k(x)$ and the degree of $k(x)$ over $k$ is denoted by $\deg(x).$ 
The symbol $\overline{x}$ always denotes the geometric point extending $x$ using the 
composition $\Spec(\overline{k})\to \Spec(k)\to X.$  If $x$ is a point of $X$ (not necessarily closed) then $\dim(x)$ denotes the dimension 
of the closure of $x.$ 
A $\bQl$-sheaf always is by definition an \'etale constructible $\bQl$-sheaf on $X$ and the associated   
derived category with bounded cohomology sheaves is denoted $\rD^b_c(X,\bQl).$
    If $x$ is a point of $X$ (not necessarily closed) 
and if $F$ is a $\bQl$-sheaf on 
$X,$ then $F_x$ denotes the restriction of $F$ to $x$ and $F_{\bar{x}}$ denotes the stalk of $F,$ viewed as a $\Gal(\overline{k(x)}/k(x))$-module.

By our assumptions on 
$k,$  the category $\rD^b_c(X,\bQl)$ is triangulated and supports Grothendieck's six operations, with internal tensor product $\otimes$  and $\Rhom,$ external product $\boxtimes,$ 
and Verdier dual 
$\vD:\rD^b_c(X,\bQl)^{\rm opp}\to \rD^b_c(X,\bQl)$
(\cite{DeligneWeil2}). For $S$ a regular scheme of dimension 
$\leq 1$ over $k$ and for a morphism of finite type $f:X\to Y$  of $S$-schemes one has the usual functors
$$ Rf_*,Rf_!: \rD^b_c(X,\bQl)\to \rD^b_c(Y,\bQl)\quad {\rm and} \quad f^*, Rf^!: \rD^b_c(Y,\bQl)\to \rD^b_c(X,\bQl)$$
with $\vD$ interchanging $ Rf_*$ and $Rf_!$ (resp.~$f^*$ and  $Rf^!$). Often one writes $f_*,f_!,$ and $f^!$ instead of  
$Rf_*,Rf_!$ and $Rf^!$ (resp.). The category of smooth (lisse) $\bQl$-sheaves on $X$ is denoted by ${\rm Lisse}(X,\bQl).$

\subsection{Remarks on perverse sheaves.} \label{subsub1}  Recall that $\rD^b_c(X,\bQl)$ contains the abelian subcategory of perverse sheaves $\Perv(X,\bQl)$  
 with respect to the autodual (middle) perversity \cite{BBD}.
 An object $K\in \rD^b_c(X,\bQl)$ is perverse if and only if the following conditions hold for any point $x\in X$ (see~\cite{BBD},~(4.0)): if $i$ denotes the 
 inclusion of $x$ into $X$ then 
 \begin{equation}\label{eqdefperv}  \cH^\nu((i^*K)_{\bar{x}})=0\textrm{ for } \nu>-\dim(x) \quad \textrm{ and }
 \quad \cH^\nu((i^!K)_{\bar{x}})=0\textrm{ for } \nu<-\dim(x).\end{equation}
 
  \begin{rem}\label{remcheckedforbark} An object  $K\in \rD^b_c(X,\bQl)$
is perverse if and only if $K|_{X\otimes\bar{k}}\in \rD^b_c(X\otimes \bar{k},\bQl)$ is perverse. 
(This is a tautology given  $i^!=\vD\circ i^*\circ \vD$ and the compatibility of 
 $\vD$ with respect  to  base change to~$\bar{k},$ cf.~\cite{KatzLaumon},~Prop.~1.1.7; \cite{BBD},~Prop.~5.1.2.) \end{rem}

 Let $j:U\hookrightarrow X$ be an open immersion with complement $i:Y\to X.$ If $K$ is a perverse sheaf on $U$ then 
there is a unique extension $j_{!*}K\in \Perv(X,\bQl)$ of $K$ to $X$ which has neither subobjects nor quotients of the form $i_*\Perv(Y,\bQl)$
(\cite{BBD}). This extension is called 
the {\it intermediate extension} or {\it middle extension}.  

Let $X$ be a smooth and geometrically connected curve over $k,$  let $j:U\hookrightarrow X$ be a dense open subscheme, and 
let $F$ be a smooth sheaf on $U.$ Then the shifted sheaf $F[1]$ (concentrated at $-1$) is a perverse sheaf on $U$ and 
the middle  extension $j_{!*}F[1]$ is a perverse sheaf which coincides with  $(j_*F)[1]$ (here we mean the usual sheaf extension $j_*F$ shifted by $1,$ see \cite{KiehlWeissauer},~Chap.~III.5).
A {\it middle extension sheaf} on $X$ ($X$ a smooth 
geometrically connected curve)
is by definition a perverse sheaf of the  form $(j_*F)[1]$ as above, cf.~\cite{Katz96},~Chap.~5.1.

\subsection{Further notions} \label{secfrobetc}
  (\cite{Laumon}) Let here $X$ denote a scheme of finite type over  $k=\FF_q$  
 and let $K\in \rD^b_c(X,\bQl).$   The geometric Frobenius morphism relative to 
 $k$ will be denoted by ${\rm Frob}_q$ or ${\rm Frob}_k.$ 

In the associated 
 Grothendieck group $\mathcal{K}(X,\bQl),$ one has an equality 
 \begin{equation}\label{eql1}  [K]=\sum_{j}(-1)^j \,[H^j(K)],\end{equation}
 with  constructible cohomology sheaves $H^j(K).$
  Recall that for any closed point 
 $x\in |X|$ and any constructible sheaf $F$ on $X,$ the stalk $F_x$ has a natural action of the geometric Frobenius 
 element $\Frob_x=\Frob_q^{\deg(k(x)/k)},$ leading to the well defined characteristic polynomial 
 $\det({1-t\cdot \Frob_x,F}).$ One defines $\tr(\Frob_x,F),$ resp. $\det(\Frob_x,F),$ 
 to be the coefficient of $-t,$ resp. $(-t)^n\,(n=\dim(F_{\bar{x}})),$ in  $\det({1-t\cdot \Frob_x,F}).$
  Using \eqref{eql1} we obtain homomorphisms of groups 
 $$ \det({1-t\cdot \Frob_x,-}):\mathcal{K}(X,\bQl)\to \bQl(t)^\times$$ $$ \tr(\Frob_x,-):\mathcal{K}(X,\bQl)\to \bQl$$ $$ \det(\Frob_x,-):\mathcal{K}(X,\bQl)\to \overline{\mathbb{\QQ}}_\ell^\times $$
by additivity/multiplicativity (cf.~\cite{Laumon}, Section~0.9). 
 This notion extends to $\rD^b_c(X,\bQl)$ by setting 
 $$ \det({1-t\cdot \Frob_x,K})=\det({1-t\cdot \Frob_x,[K]}).$$

 Let $X$ be a curve and let $F$ be a smooth $\bQl$-sheaf on a dense open subset $j:U\hookrightarrow X.$ If $x\in |X|$ then $X_{(x)}$ (resp. $X_{(\overline{x})}$) 
denotes the Henselization 
of $X$ with respect to $x$ (resp. $\bar{x}$) and $\eta_x$ (resp. $\overline{\eta}_x$) denotes the generic point of $X_{(x)}$ (resp. $X_{(\bar{x})}$), cf.~\cite{DeligneWeil2}.   
One defines the {\it generic rank} $r(F)=r(j_*F)$ of $F$ as $\rk(F_{\eta_x})\, (x\in X)$  and extends this notion to $K\in \rD^b_c(X)$
by additivity, cf.~\cite{Laumon},~2.2.1.

\subsection{Artin-Schreier and Kummer sheaves}\label{subsub2} Recall the construction of Artin-Schreier and Kummer sheaves: 
Let $k$ be the finite field $\FF_q$ and let $G$ be a commutative connected algebraic group  of finite type over $k.$ 

\begin{rem}\label{rempullback}
If $f:X\to G$ is a morphism of schemes, and if $\cL$ is a sheaf on $G$ then we set $\L(f)=f^*\L.$ 
Sometimes we simply $\L$ instead of 
 $\L(f),$ especially if $f$ is an obvious change of base. 
\end{rem}

The Lang isogeny of $G$ is the extension of $G$ by $G(k)$ $$ 1\to G(k)\to G\stackrel{L}{\to}G\to 1$$ 
where $L(x)={\rm Frob}_q(x)\cdot x^{-1}.$  Hence $L$ exhibits $G$ as a $G(k)$-torsor over itself, the {\it Lang torsor}. 
To a character $\chi:G(k) \to \overline{\mathbb{\QQ}}_\ell^\times$ one then associates a smooth rank-one sheaf $\cL_\chi$  on $G$ 
by pushing out 
the Lang torsor 
by $\chi^{-1}: G(k) \to \overline{\mathbb{\QQ}}_\ell^\times.$  

If $G=\GG_{m,k},$ then $\L_\chi$ is called a {\it Kummer sheaf}. If $G=\AA^1_k,$ then
a nontrivial character  $\AA^1(k)\to \overline{\mathbb{\QQ}}_\ell^\times$ is usually denoted by $\psi$ and  the resulting sheaf 
 $\L_\psi$ is called 
an {\it Artin-Schreier} sheaf.

Consider 
the multiplication map $$x\cdot x':\AA^1\times_k \AA^1\to \AA^1,\, (x,x')\mapsto x\cdot x'.$$ Then, for a closed point $s$ of $\AA^1,$  the restriction of $\cL_\psi(x\cdot x')$ to $s\times_k \AA^1$ is denoted by $\cL_\psi(s\cdot x').$ 

%For $G=\GG_m$ this construction carries over  to the case where $k$ is algebraically closed, cf.~\cite{Laumon}, Rem.~1.1.3.7,
%so we obtain a Kummer sheaf $\cL_\chi$ on $\GG_{m,\overline{k}}$ associated to a character $\chi: \GG_m(k_0)\to \overline{\mathbb{\QQ}}_\ell^\times,$ where $k_0$ is a finite
%subfield of $\overline{k}$. 
If $k$ is a field of odd order then 
the unique quadratic character  $\GG_m(k)=k^\times\to \overline{\mathbb{\QQ}}_\ell^\times $ is denoted~$-\1.$ 
The trivial character $\GG_m(k)=k^\times\to \overline{\mathbb{\QQ}}_\ell^\times $ is denoted~$\1.$ 

\begin{rem}  Recall from \cite{Laumon},~(1.1.3.7), an alternative construction of Kummer sheaves:
 let $k$ be momentarily allowed to be any field which contains a primitive $N$-th root of unity. Consider the exact sequence
 $$ 1\to \mu_N(k)\to \GG_{m,k}\stackrel{[N]}{\longrightarrow}\GG_{m,k}\to 1,$$
 where $[N]$ denotes the $N$-th power morphism and $\mu_N(k)$ is the group of 
 $N$-th roots of unity in $k.$ 
 
 This is a Galois cover with Galois group $\mu_N(k)$ 
 and pushing out the resulting $\mu_N(k)$-torsor  via $\chi^{-1}$ (where  $\chi:\mu_N(k)\to \bQl^\times$ is a character) one 
 obtains a Kummer sheaf $\mathcal{K}_\chi$ on $\GG_{m,k}.$ Then the  sheaf $\mathcal{K}_\chi$
  coincides with the above Kummer sheaves 
 $\cL_\chi$ if $k=\FF_q$  and $N=q-1$ (\cite{Laumon},~(1.1.3.7)). 
\end{rem}

\section{Convolution in characteristic $p$}\label{secvondefetale2} 
\subsection{Basic definitions.}
In this section $k$ denotes either a finite field of characteristic $p\neq \ell$ or the algebraic
closure of such a field.  
Let us recall the definitions and basic results of \cite{Katz96},~Section~2.5. 
For $G$ a smooth $k$-group, denote the multiplication map by $\pi:G\times G\to G.$ Let $K$ and $L$ be two objects of $\rD^b_c(G,\bQl)$
and let $K\boxtimes L$ denote the external tensor product of $K$ and $L$ on $G\times G$ with respect to the two natural projections.  
Then one may form the {\it $!$-convolution} $$K*_!L:=R\pi_!(K\boxtimes L)$$ as well as the  {\it $*$-convolution}   $$K*_*L:=R\pi_*(K\boxtimes L) $$
with duality interchanging both types of convolution. 
Under the shearing transformation  $$\sigma:\AA^2_{x,y}\to \AA^2_{x,t},\, (x,y)\mapsto (x,t=x+y),$$ the above convolutions can be written as 
$$K*_!L:=R\pr_{2!}(K\boxtimes L),\quad K*_*L:=R\pr_{2*}(K\boxtimes L), $$ where the external tensor product is now
formed with respect to the first projection $\pr_1:\AA^1_x\times \AA^1_t\to \AA^1_x$ and the difference map
$$ \delta:\AA^1_x\times \AA^1_t\to \AA^1_y,\,(x,t)\mapsto y=t-x.$$ 

An object $K$ of $\Perv(G,\bQl)$ has property $\cP$ by definition if 
for any perverse sheaf $L\in \Perv(G,\bQl)$ the convolutions $L*_!K$ as well as $L*_*K$ are again perverse. 
If either $K$ or $L$ has the property 
$\cP$ then one can define the {\it middle 
convolution} of $K$ and $L$ as the image of $L*_!K$ in $L*_*K$ under the natural forget supports map 
\begin{equation} \label{eqdefmiddleconv} L*_\middle K:=\im(L*_!K\to L*_*K).\end{equation}

%It turns out that  the middle convolution on the affine line admits a concrete description in terms of 
%a variation of ``parabolic" cohomology groups (cf.~\cite{Katz96},~Cor.~2.8.5). We need two preparatory results:

\begin{lem}\label{lemP}  Let $K$ be a perverse sheaf on $\AA^1_k$ which is geometrically irreducible 
and not geometrically translation invariant. Then $K$ has the property $\cP.$ 
\end{lem}
\proof 
This follows from \cite{Katz96},~Cor.~2.6.10 together with Rem.~\ref{remcheckedforbark}. 
\Endproof

Using the previous result we obtain for each Kummer sheaf $\L_\chi,$  associated to a {\it nontrivial} character $\chi,$  a functor
\begin{equation}\label{eqmcchidef}  \MC_\chi: \Perv(\AA^1_k,\bQl)\to \Perv(\AA^1_k,\bQl),\, K\mapsto K*_\middle L_\chi,\end{equation}
with $L_\chi=j_*\L_\chi[1],$ where $j:\GG_m\to \AA^1$ denotes the natural inclusion. \\

Let now $S$ be any $k$-variety, let $\overline{f}:X\to S$ be proper, let $j:U\to X$ be an affine open immersion over $S,$ let $D=X\setminus U,$  and 
suppose that $\overline{f}|_{D}: D\to S$ is affine. 
Suppose that $K$ is an object in $\Perv(U,\bQl)$ such that both $Rf_!K$ and $Rf_*K$ are perverse. Then Prop.~2.7.2 of \cite{Katz96} states 
that $R\overline{f}_*(j_{!*}K)$ is again perverse
and that \begin{equation}\label{eqmiddleext}
     R\overline{f}_*(j_{!*}K) =\im(Rf_!K\to Rf_*K).
    \end{equation}
Let us take  
$$S=\AA^1_t,\quad X=\PP^1_x\times\AA^1_t,\quad U=\AA^2_{x,t},$$ and let 
$f=\pr_2:\AA^2_{x,t}\to \AA^1_t,$ and $\overline{f}=\overline{\pr}_2:\PP^1_x\times\AA^1_t\to \AA^1_t.$  
Then Eq.~\eqref{eqmiddleext} implies:

\begin{lem}\label{lemmiddle} Let $K\in \Perv(\AA^1,\bQl)$ have property $\cP$ and let 
$L\in \Perv(\AA^1,\bQl ).$ Then $K*_\middle L$ is a perverse sheaf with 
\begin{equation}\label{eqmiddleext2}
 K*_\middle L=R\overline{\pr}_{2*}(j_{!*}(K\boxtimes L))\quad {\rm with} \quad K\boxtimes L=\pr_1^*K\otimes \delta^*L.
\end{equation}
 \end{lem}

\begin{prop}\label{rembasicprops}(Katz) 
The middle convolution has the following properties:
\begin{enumerate}
\item  Let $K\in \Perv(\AA^1_x,\bQl)$ and $L\in \Perv(\AA^1_y,\bQl)$ 
be irreducible middle extensions which 
are not geometrically translation invariant.If $K$ and $L$ are tame at $\infty$ then there is a short exact sequence of perverse sheaves on $\AA^1_t$
$$ 0\to H\to K*_!L\to K*_\middle L \to 0,$$
 where $H$ is the constant sheaf $\pr_{2*}(j_{*}(K\boxtimes L)_{\infty \times \AA^1_t})$ on $\AA^1_t$ (\cite{Katz96}, 2.9.4). 
 \item If $F,K,L\in \rD^b_c(\AA^1,\bQl)$ have all property $\cP$ then 
 $$ F*_\middle(K*_\middle L)=(F*_\middle K)*_\middle L,$$ 
 cf.~\cite{Katz96},~2.6.5.
 \item For each nontrivial Kummer sheaf $\cL_\chi$  and for each $K\in \Perv(\AA^1,\bQl)$ having the property $\cP,$ the following holds:
 $$  \MC_{\chi^{-1}}(\MC_\chi(K))=K(-1).$$ 
 This follows from~(i) using $L_{\chi^{-1}}*_\middle L_\chi=\delta_0(-1)$ with $L_\chi=j_*\cL_\chi[1]$ and with $\delta_0$ denoting the trivial 
 sheaf supported at $0,$ cf.~\cite{Katz96},~Thm.~2.9.7. 
% \item If $K,L\in \Perv(\AA^1,\bQl)$ have both the property $\cP$ then $K*_\middle L$ has again the property $\cP.$ 
\item Formation of $K*_\middle L$ is compatible with arbitrary change of base (see \cite{Katz96}~4.3.8--4.3.11). 
\end{enumerate}

\end{prop}\Endproof

\subsection{Fourier transformation and convolution.}\label{secintrofourier}  In this section, we fix a finite field $k=\FF_q\, (q=p^m)$ and an  additive $\overline{\mathbb{\QQ}}_\ell^\times$-character 
$\psi$ of $\AA^1(\FF_p),$ inducing for all $k\in \NN$ an additive character $\psi_{\FF_{q^k}}=\psi=\psi\circ \tr^{\FF_{q^k}}_{\FF_q}.$

By the discussion in Section~\ref{subsub2} we have the associated Artin-Schreier sheaf $\cL_\psi$ on $\AA^1_k.$ 
Let $\AA=\Spec(k[x])$ and $\AA'=\Spec(k[x'])$ be two copies of the affine line and let 
$$ x\cdot x':\AA\times \AA'\To \GG_{a,k},\quad (x,x')\mapsto x\cdot x' .$$ 
The two projections of $\AA\times \AA'$ to $\AA$ and $\AA'$ are denoted $\pr$ and $\pr',$ respectively.    
Following Deligne and Laumon \cite{Laumon}, we can form the Fourier transform as follows:
$$\cF_\psi=\cF: \rD^b_c(\AA,\bQl)\To \rD^b_c(\AA',\bQl),\quad K\Mapsto R\pr'_!\left(\pr^*K\otimes \cL_\psi(x\cdot x')\right)[1].$$ 
By exchanging the roles of $\AA$ and $\AA',$ one obtains the Fourier transform
$$\cF_\psi'=\cF': \rD^b_c(\AA',\bQl)\To \rD^b_c(\AA,\bQl),\quad K\Mapsto R\pr_!\left(\pr^{'*}K\otimes \cL_\psi(x\cdot x')\right)[1].$$ 
Consider the automorphism $a:\AA\to \AA, x\mapsto -x.$ 
By \cite{Laumon},~Cor.~1.2.2.3 and Thm.~1.3.2.3, the Fourier transform is an equivalence of triangulated categories 
$\rD^b_c(\AA,\bQl)\to \rD^b_c(\AA',\bQl)$ and $\Perv(\AA,\bQl)\to \Perv(\AA',\bQl)$ with quasi-inverse  $a^* \cF'(-)(1).$ 
Especially, it maps simple objects to simple objects.
\begin{defn}{\rm 
Let ${\rm Fourier}(\AA,\bQl)\subset \Perv(\AA,\bQl)$ and ${\rm Fourier}(\AA',\bQl)\subset \Perv(\AA',\bQl)$ be the categories of simple middle extension sheaves  on $\AA_k$ and $\AA'_k$ (resp.) which are not geometrically isomorphic to 
a translated Artin-Schreier sheaf $\cL_\psi(s\cdot x)$ with $s\in \bar{k}$ (cf.~\cite{Laumon},~(1.4.2)). We call the objects in ${\rm Fourier}(\AA,\bQl)$ {\it irreducible Fourier sheaves}.}
\end{defn}

             In \cite{Katz90},~(7.3.6), the sheaves $\cH^{-1}(K)$ with $K\in {\rm Fourier}(\AA,\bQl)$ 
             are called {\it irreducible Fourier sheaves}, justifying the nomenclature (up to a shift).  
By 
Thm.~1.4.2.1 and Thm.~1.4.3.2 in \cite{Laumon}, the following holds:

\begin{prop}\label{propkatz1}(Deligne, Laumon) 
\begin{enumerate}
\item 
The functor $\cF$  induces a categorial  equivalence from  ${\rm Fourier}(\AA,\bQl)$ to ${\rm Fourier}(\AA',\bQl).$ 
\item If $H=V\otimes \cL_\psi(s\cdot x),\, (s\in |\AA^1|)$ with $V$ constant, then  $\cF_\psi(H)$ is the punctual sheaf $ V_s$ supported 
at $s.$  
\item Let $\chi$ be a nontrivial character of $\GG_{m}(k)$  
and let $j:\GG_m\to \AA,$ resp. $j':\GG_m'\to \AA'$ denote the canonical inclusions. Then 
$$ \cF(j_*\cL_\chi[1])=j'_*\cL_{\chi^{-1}}[1]\otimes G(\chi,\psi),$$
where $G(\chi,\psi)$ is the geometrically constant sheaf   on $\AA'$ on which the Frobenius acts via the Gau{\ss} sum
$$ g(\chi,\psi)=-\sum_{x\in k^\times} \chi(x)\psi(x)$$
(as a $\Frob_q$-module, $G(\chi,\psi)=H^1_c(\AA_{\overline{k}}\setminus 0, \cL_\chi\otimes(\cL_\psi|_{\AA^1\setminus 0}))$). 
\end{enumerate}
\end{prop}

\begin{rem}  An irreducible perverse sheaf $K\in \Perv(\AA,\bQl)$ has the property $\cP$ if  and only if $\cF(K)$ is a middle extension
(cf.~\cite{Katz96},~2.10.3). Note that the trivial rank-one sheaf $\bQl$ can be viewed as $\cL_\psi(0\cdot x').$ It follows hence 
from by Lem.~\ref{lemP} that  any object in ${\rm Fourier}(\AA,\bQl)$ and in ${\rm Fourier}(\AA',\bQl)$ has the property~$\cP.$
\end{rem}

The relation of the Fourier transform to the convolution is expressed as follows (\cite{Laumon},~Prop.~1.2.2.7):
\begin{equation}\label{eqfourier}
 \cF(K_1*_!K_2)=(\cF(K_1)\otimes \cF(K_2))[-1]\,.
 %\quad {\rm and} \quad \cF(K_1\otimes K_2)[-1]=\cF(K_1)*_!\cF(K_2).
\end{equation}
Applying Fourier inversion  yields
\begin{equation}\label{eqfourier2}
 K_1*_!K_2=a^*\cF'(\cF(K_1)\otimes \cF(K_2))[-1](1)\,.
\end{equation} 

\begin{prop}\label{proprelfourier} Let $K,L\in \Perv(\AA_k,\bQl)$ be tame   middle extensions in ${\rm Fourier}(\AA,\bQl).$ Suppose that for 
$j:\GG_m\hookrightarrow \AA'$ the inclusion one has 
$$\cF(K)=j_*F[1]\in {\rm Fourier}(\AA',\bQl)\quad \textrm{and}\quad  \cF(L)=j_*G[1]\in {\rm Fourier}(\AA',\bQl)$$
for smooth sheaves $F,G$ on $\GG_m.$  Then the following holds:
\begin{enumerate}
 \item 
$$ \cF(K*_\middle L)=j_*(F\otimes G)[1]\, .$$
\item if  $K$ is not a translate of $j_*\cL_{\chi^{-1}}[1],$ then $\cF(K*_\middle j_*\cL_{\chi^{-1}}[1] )$ is an object in ${\rm Fourier}(\AA',\bQl).$ 
\end{enumerate}
\end{prop}

\proof This is \cite{Katz96},~2.10.8 (the proof generalizes to our context).
%It follows from Thm.~\ref{thmKatz1}(ii) that there is a short exact sequence of perverse sheaves
%$$ 0\to H\to K\ast_!L\to K\ast_\middle L\to 0$$ with $H$ a constant sheaf shifted by $1.$  
%The exactness of Fourier transform together with  Prop.~\ref{propkatz1}(ii) and \eqref{eqfourier} give an exact sequence  
%$$ 0\to \textrm{punctual sheaf, supported at $0$}\to \cF(K*_!L)=(\cF(K)\otimes \cF(L))[-1]\to \cF(K*_\middle L)\to 0.$$ 
%Hence, over $\GG_m,$ the restriction of the above sequence gives 
%$$ j^*(\cF(K)\otimes \cF(L))[-1]=(F\otimes G)[1]=j^*\cF(K*_\middle L).$$ It follows from \cite{Katz96}, Cor. 2.6.17, and 
%from Rem.~\ref{remcheckedforbark} that $K*_\middle L$ has again the property $\cP$ 
%which implies that $\cF(K*_\middle L)$ is a middle extension by the remark following the definition of ${\rm Fourier}(\AA',\bQl).$ 
%Hence we obtain 
%$$ \cF(K*_\middle L)=j_*((F\otimes G)[1]),$$ proving the first claim. 
The second claim holds since, under the given assumptions
on $K$ and $L,$ the sheaf $j_*((F\otimes G)[1])$ is irreducible and not an Artin-Schreier sheaf. 
\Endproof

\begin{cor}\label{corfourier} Under the assumptions of Prop.~\ref{proprelfourier}:
 $$ K*_\middle L =a^*\cF'\left(j_*(F\otimes G\right)[1])(1)\,.$$ Moreover, if 
 $L=j_*\cL_\chi[1]$ and if $L$ is not a translate of  $j_*\cL_{\chi^{-1}}[1],$ then $ K*_\middle L\in {\rm Fourier}(\AA,\bQl).$ 
\end{cor}

\proof This follows from Fourier inversion and from Prop.~\ref{proprelfourier}~(ii). \Endproof

\section{Local Fourier transform and local monodromy of the middle convolution.}
\subsection{Local Fourier transform.}  As before, we fix a finite field $k=\FF_q\, (q=p^m),$
a prime $\ell\neq p,$  and an  additive $\overline{\mathbb{\QQ}}_\ell^\times$-character 
$\psi$ of $\AA^1(\FF_q),$ see Section~\ref{secintrofourier}.  
In the following we summarize  Laumon's construction of the local Fourier transform \cite{Laumon} and 
the stationary phase decomposition:\\

 Recall the notion of a {\it henselian trait} $T=(T,\eta,t)$  (see~\cite{DeligneWeil2}, \cite{Laumon}): it is the spectrum of a henselian discrete valuation ring $R$ having   generic point $\eta$ and closed point $t=\Spec(k).$ The usual geometric  point over $t$  is denoted $\bar{t}=\Spec(\bar{k})$ 
 and the generic point of the strict henselization $T_{\bar{t}}$ of $T$ 
 is  denoted $\bar{\eta}$ (\cite{DeligneWeil2},~(0.6)).   Note that any henselian trait $T$ (resp. $R$) has a uniformizer, usually denoted~$\pi.$
 
 Below we  will consider two henselian traits
 $T=(T,\eta,t)$ and $T'=(T',\eta',t')$   in equiconstant characteristic $p$ with given 
 uniformizers $\pi,$ resp. $\pi',$ having $k$ as residue field. The fundamental groups $\pi_1(\eta,\overline{\eta})$ and  $\pi_1(\eta',\overline{\eta}')$ are denoted $G$ and $G',$ respectively. The inertia subgroups of $G,G'$ are denoted $I,I'$ (resp.). \\

 For $X$ a variety
 over $k$ and $x$ a closed point of $X$ one has the henselian trait 
 $X_{(x)}=(X_{(x)},\eta_x,x)$ given by 
 the spectrum of the henselization $R_x$ of the local ring $\mathcal{O}_{X,x},$ its generic point $\eta_x$ and the point~$x.$ We define $G_x=\pi_1(\eta_x,\overline{\eta}_x)$ 
 and  $I_x\leq G_x$ denotes the inertia subgroup of $G_x.$  \\

 The category of smooth $\bQl$-sheaves on $\eta$ (which will
  be identified  with the category 
 of continuous $\bQl$-representations of finite rank of $G$ via the fibre functor 
 in $\bar{\eta},$  cf. \cite{Laumon},~Rem.~2.1.2.1) 
 is denoted $\cG.$ Similarly 
 we define the category $\cG'$ of smooth sheaves on $\eta'.$
 For $V\in \ob \cG,$ denote by $ V_!$ the extension by zero to $T,$ similarly for $V'\in \cG'.$   
 The subcategory of $\cG,$ resp. $\cG',$ formed by objects whose 
 inertial slopes are in $[0,1[$ are denoted $\cG_{[0,1[},$ resp.
 $\cG'_{[0,1[},$ cf.~\cite{Laumon}, Section~2.1. 
 Recall that an object of $\cG$ is tamely ramified if and only if 
 it is pure of slope $0$ (loc.cit., 2.1.4).
  If $V$ 
 (resp. $V'$) is an object of $\cG$
 (resp. of $\cG'$) then its extension by zero to $T$ (resp. $T'$) is denoted
 $ V_!$ (resp. $V'_!$).\\
 
 One has the $\bQl$-sheaves ${\cL}_\psi(\pi/\pi'),\,{\cL}_\psi(\pi'/\pi)$ and ${\cL}_\psi(1/\pi\pi')$
 on $T\times_k\eta',\, \eta\times_kT'$ and $\eta\times \eta'$ (resp.) and the respective extensions by zero to 
 $T\times_k T'$ are denoted   $\bar{\cL}_\psi(\pi/\pi'),\,\bar{\cL}_\psi(\pi'/\pi)$ and $\bar{\cL}_\psi(1/\pi\pi').$
  For any $V\in \ob \cG$ one may form the vanishing cycles 
 $$ R\Phi_{\eta'}(\pr^* V_!\otimes \bar{\cL}_\psi(\pi/\pi')),\quad R\Phi_{\eta'}(\pr^* V_!\otimes \bar{\cL}_\psi(\pi'/\pi)),
 \quad R\Phi_{\eta'}(\pr^* V_!\otimes \bar{\cL}_\psi(1/\pi'\pi))\quad$$ as objects in  $ \rD^b_c(T\times_k \eta',\bQl)$
 with respect to $\pr':T\times_k T'\to T'$ ($\pr:T\times_k T'\to T$ denoting the first projection), see \cite{DeligneVanishing},~(2.1.1). \\
 
 This leads to three functors, called 
 {\it local Fourier transforms},  
 $$\cF^{(0,\infty')},\cF^{(\infty,0')},\cF^{(\infty,\infty')}: \cG\to \cG',$$ defined by 
 $$ \cF^{(0,\infty')}(V)=R^1\Phi_{\bar{\eta}'}(\pr^* V_!\otimes \bar{\cL}_\psi(\pi/\pi'))_{(\bar{t},\bar{t}')}\,,$$ 
 $$ \cF^{(\infty,0')}(V)=R^1\Phi_{\bar{\eta}'}(\pr^* V_!\otimes \bar{\cL}_\psi(\pi'/\pi))_{(\bar{t},\bar{t}')}\,,$$ 
 $$ \cF^{(\infty,\infty')}(V)=R^1\Phi_{\bar{\eta}'}(\pr^* V_!\otimes \bar{\cL}_\psi(1/\pi\pi'))_{(\bar{t},\bar{t}')}\,,$$ 
 see~\cite{Laumon},~2.4.2.3. By interchanging the roles of $T$ and $T'$ one obtains the functors $$\cF^{(0',\infty)},\cF^{(\infty',0)},\cF^{(\infty',\infty)}: \cG'\to \cG.$$

 \begin{rem}\label{remvariousT}
 Note that we have neither fixed $T$ nor $T'$ so that the local 
 Fourier transform may be formed with respect to any pair of henselian traits in equiconstant characteristic $p$ 
  having  $k$ as residue field.
  \end{rem}
  
   We will need the following properties of the local Fourier transform below:
 
\begin{thm}\label{thmlaumon} (Laumon) 
\begin{enumerate}
\item The three local Fourier transforms are exact functors. Moreover,
 $\cF^{(0,\infty')}:\cG\to \cG'_{[0,1[}$ is an equivalence of categories quasi-inverse
to $a^*\cF^{(\infty',0)}(-)(1),$ where $a:T\to T$ is the automorphism defined by 
$\pi\mapsto -\pi$ and $(1)$ denotes a Tate-twist. 
\item If $W$ denotes an unramified  $G$-module, then 
$$ \cF^{(0,\infty')}(W)=W,\quad \cF^{(\infty,0')}(W)=W(-1),\quad \cF^{(\infty,\infty')}(W)=0.$$
\item For a non-trivial Kummer sheaf $\cK_\chi$ on $\GG_m=
\Spec(k[u,u^{-1}]),$ denote $  V_\chi,$ resp.
$V'_\chi$ the $G$-module $\cK_\chi(\pi)$ (resp. the $G'$-module $\cK_\chi(\pi')$) on 
$T$ (resp. $T'$), where $\pi:\eta\to \GG_{m}$ (resp; $\pi':\eta'\to \GG_m$) is the morphism which 
maps $\pi$ to $u$ (resp. $\pi'$ to $u$). Then
$$ \cF^{(0,\infty')}( V_\chi)=V'_\chi\otimes G(\chi,\psi),$$
$$ \cF^{(\infty,0')}( V_\chi)=V'_\chi \otimes G(\chi^{-1},\psi),$$
where $G(\chi,\psi)$ denotes the unramified $G$-module $H^1_c(\GG_{m,\overline{k}},\cK_\chi\otimes 
\cL_\psi)$ whose Frobenius trace is the Gau{\ss} sum 
$$\tr(\Frob_k,G(\chi,\psi))=g(\chi,\psi)=-\sum_{a\in k^\times}\chi(a)\psi_k(a).$$
\item If the restriction of the representation $V$ to the inertia subgroup 
$I$ is unipotent indecomposable (resp. tame),
then $\cF^{(0,\infty')}(V),$ resp. $\cF^{(\infty,0')}(V),$ is   unipotent and indecomposable (resp. tame) of the same rank.
\item If $W$ is an unramified $G$-module and if $W'$ denotes 
the unique $G'$-module corresponding to the same ${\rm Gal}(\bar{k}/k)$-module then
$$ \cF^{(0,\infty')}( V_\chi\otimes W)=\cF^{(0,\infty')}( V_\chi)\otimes W'.$$ 
\item Let  $T_1=T\otimes_k k_1$ with $k_1$ a finite extension of $k,$ let $\eta_1$ denote the generic point of $T_1$ and let 
$G_1=\Gal(\overline{\eta}_1/\eta_1).$ Let $f:T_1\to T$ denote the \'etale map given by the canonical projection. 
If $V$ is a tamely ramified irreducible $G$-module of the form 
$V={\rm Ind}_{G_1}^G( V_1),$ for $ V_1$ a rank-$1$ module of $G_1$ then the following holds:
$$ \cF^{(0,\infty')}(V)={\rm Ind}_{G_1}^G(\cF^{(0_1,\infty_1')}( V_1)).$$
\end{enumerate}
\end{thm}

\proof The assertions (i)--(iii) are contained in \cite{Laumon}~Thm.~2.4.3, Prop.~2.5.3.1. Assertion~(iv) is proven in 
\cite{KatzTraveauxdeLaumon}, Lemma~5. Assertion~(v) can be found in \cite{Laumon},~(3.5.3.1). Assertion (vi)
follows from proper base change (\cite{Laumon},~(2.5.2), (3.5.3.1)).
\Endproof

 \begin{rem} \label{remgaussformula} {\it The following formula for Gau{\ss} sums will be used below:
  \begin{equation}\label{chichiinv} g(\chi,\psi)g(\chi^{-1},\psi)=\chi(-1)\cdot q\, .\end{equation}}
 \end{rem}

The following result is  one version of Laumon's principle of stationary phase \cite{Laumon}.  We use the formulation of Katz in  \cite{Katz90},~Cor.~7.4.2. Although in loc.~cit. the result is stated 
for $k=\bar{\FF}_q,$ this is not essential (as remarked in \cite{Katz90}, beginning of Section~7.4, this is only for notational reasons) and the transition to $k=\FF_q$ can be made using the arguments of
\cite{Laumon}, preuve de (3.4.2). 
See Rem.~\ref{remvariousT} for an explanation of the modules $\cF^{(0,\infty')}(F_{\overline{\eta}_s}/F_{\overline{\eta}_s}^{I_s^t})$ below (using the uniformizer $\pi=x-s$).

\begin{thm}\label{corLaumonKatz}(Laumon)  Let $K=j_*F[1]\in {\rm Fourier}(\AA^1_k,\bQl)$ be a middle extension of a  smooth  sheaf $F$ on 
$ U=\AA^1\setminus S\stackrel{j}{\hookrightarrow}\AA^1$ which is tamely ramified  at $S\cup \infty$ and let $F'[1]=\cF(K).$  
Then  there is an isomorphism of $G_{\infty'}$-modules
\begin{equation}\label{eqstationary2} F'_{\overline{\eta}_\infty}\simeq \bigoplus_{s\in S}{\rm Ind}_{G_{s\times_k\infty'}}^{G_{\infty'}}
(\cF^{(0,\infty')}(F_{\overline{\eta}_s}/F_{\overline{\eta}_s}^{I_s^t})\otimes \bar{\cL}_\psi(s\cdot x')_{\overline{\eta}_{\infty'}}) \,.\end{equation}  \end{thm}\Endproof

\begin{cor}\label{corstationaryback} Let $K\in {\rm Fourier}(\AA^1_k,\bQl)$ be a middle extension of a  smooth sheaf $F$ on 
$ U\stackrel{j}{\hookrightarrow}\AA^1_k,$ tamely ramified at $S\cup \infty.$  Let $K':=\cF(K)=F'[1].$  Write  the stationary phase decomposition in 
\eqref{eqstationary2} as
\begin{equation}\label{eqstationary3} F'_{\overline{\eta}_\infty}\simeq \bigoplus_{s\in S}{\rm Ind}_{G_{s\times_k\infty'}}^{G_{\infty'}}
(V'_s \otimes \bar{\cL}_\psi(s\cdot x')_{\overline{\eta}_{\infty'}}) \,.\end{equation} 
For each $s\in S$ there is an isomorphism of $G_s$-modules
$$ F_{\overline{\eta}_s}/ F_{\overline{\eta}_s}^{I_s}\simeq a^*\cF^{(\infty',0)}(V'_s)(1).$$
%where the local Fourier transform $\cF^{(\infty',0)}$ is formed with respect to $s\times D'_{\bar{\eta}_{\infty}}$ and $\AA_{(s)}.$  
 \end{cor}
 
 \proof Using Thm.~\ref{thmlaumon}(vi) one reduces to the case where $k(s)=k.$ 
 Suppose first that $s=0.$ Then \cite{Katz90}, Cor.~7.4.3.1, states that 
 $$ \cF(F')_{\bar{\eta}_0}/\cF(F')_{\bar{\eta}_0}^{I_0}\simeq \cF^{(\infty',0)}(F'_{\overline{\eta}_\infty}).$$
 Since in \eqref{eqstationary2} the summand belonging to $s=0$ is uniquely determined by 
 its slope being equal to $0$ (see \cite{Katz90}, Cor. 7.4.1.1) and since 
 the other summands in have slope 
 equal to $1$ (loc.cit.) it follows from $\cF^{(\infty',0)}(W)=0\,\forall W\in \cG_{[1,\infty[}$
 (see \cite{Laumon}, Thm.~(2.4.3)(ii)b)) that 
 $$ \cF^{(\infty',0)}(F'_{\overline{\eta}_\infty})=\cF^{(\infty',0)}(V'_{0}).$$ The claim follows 
 now from Fourier inversion (see 
 Prop.~\ref{propkatz1}(i) and Thm.~\ref{thmlaumon}(i)) for $s=0.$   If $s\neq 0$ then 
 (as in \cite{Katz90}, 7.4.1) one 
 can use the formula  
 $$ \cF(F)\simeq \cF({\rm Add}(s)^*(F))\otimes \cL_{\psi(sx')},$$ 
 where ${\rm Add}(s): x\mapsto x+s,$ in order to reduce to the case $s=0.$ 
 \Endproof
%
%We remark that in the last two results one may relax the assumption from tameness at $S\cup \infty$ to tameness at $\infty.$ 

 \subsection{Local monodromy of the middle convolution with Kummer sheaves.}\label{secThom}
 Let $T$ be a henselian trait with residue field $k=\FF_q,$ with uniformizer $\pi,$  
generic point $\eta,$ closed point $t,$ and with fraction field $K_t.$ 
 Then the tame quotient $G^t$ of the fundamental group $G=\pi_1(\eta, \overline{\eta})$
is a semidirect product of the procyclic tame inertia group \begin{equation}\label{eqzl1}I^t\simeq 
  \hat{\ZZ}(1)(\overline{k}):=\lim_{\stackrel{\longleftarrow}{N\geq 1, (N,p)=1}}\mu_N(\bar{k}),\end{equation} and the absolute Galois group 
  ${\rm Gal}(\bar{k}/k)$ of the residue field of $T,$
  cf.~\cite{Laumon}, Section~2.1.  \\
%  \begin{rem}\label{remkummersheaves}{\it  Let $\cL_\chi$ be  a 
%Kummer sheaf on $\GG_{m,k}\subset \AA^1_k$ 
%as defined in   Section~\ref{subsub2}. Let for the moment $T=\AA^1_{(0)}$ and let us 
%denote the restriction of $\cL_\chi$ to the generic point $\eta$ of $T$ again by $\cL_\chi.$ Then $\cL_\chi$ 
% corresponds to a  character $\rho_\chi$ of the abelianization of $G^t=\hat{\ZZ}(1)(\bar{k})\rtimes \Gal(\bar{k}/k)$ (the abelianization being 
%  isomorphic to the direct product of 
%$k^\times$  
%and $\Gal(\bar{k}/k)$).  
%%Then  the very construction 
%%of $\cL_\chi$ implies that $\rho_\chi(\Frob_k)=1.$ 
% }  \end{rem}
%  

  For $l\in \NN_{>1}$ let $k_l=\FF_{q^l},$ let $T_l:=T\times_k   k_l,$ with $T_l$ having 
residue field $k_l$ and generic point $\eta_l.$  Let 
$G_l^t$ denote the tame quotient of $\pi_1(\eta_l,\bar{\eta}_l)\simeq \Gal(\bar{\eta}/\eta_l),$ semidirect product of $I^t$ and $\Gal(\bar{k}/{k_l}).$ We view $G_l^t$ as a subgroup of $G^t$ in the obvious way.\\

 Moreover, for $s$ a closed point of $\AA^1$ with henselian trait 
$(T_s=\AA^1_{(s)}, \eta_s, s)$ we set $G_s:=\pi_1(\eta_s,\bar{\eta}_s)$ and 
$G_s^t$ denotes the tame quotient of $G_s.$ Note that since $s\simeq \Spec(\FF_{q^l})$
for some $l$ there is an isomorphism between  $G_s^t$ and this $G_l^t.$ \\

  By the theorem of Krull, Remak and Schmidt,
any   $G^t$-module $W$
  decomposes 
  into a direct sum of indecomposable summands $ V_1\oplus \cdots \oplus  V_k,$ unique up to renumeration.  
 In the following we suppose that each indecomposable summand $V$ of $W$ is of the form
   \begin{equation}\label{eqNFiltr2} V=\cJ_{n}\otimes \Ind_{G_l^t}^{G^t}( V_{\chi'}\otimes F)\end{equation}
  with 
  \begin{enumerate}
  \item $ V_{\chi'}$ denoting the rank one $G_l^t$-module of  the  Kummer sheaf belonging to  $\chi':k_l^\times \to \overline{\mathbb{\QQ}}_\ell^\times$ (as in Thm~\ref{thmlaumon}(iii)), 
  \item  $F$ an unramified $G^t_l$-module
  of rank~$1,$ 
  \item  $\cJ_{n}$ an
  indecomposable $G^t$-module of rank~$n\geq 1$ on which the group $I^t$ acts unipotently
   such that the following holds:
   %the $I^t$-eigenspace is isomorphic as a $\Frob_k$-module to 
%   $\bQl(1-n).$ 
%     Note $\cJ_n$ is not uniquely determined by the latter two conditions but, due to the 
%  theory of Jordan normal forms, the associated monodromy filtrations, as defined below,  behave similarly (\cite{DeligneWeil2},~(1.6.7.1)).
%  So $\cJ_n$ stands for a class of representations rather than a unique representation. 
%  
%  
%  Consider the Frobenius-weight filtration $M$ on a $G^t$-module $V,$ compatible with the operation 
%  of the associated Weil group
%  (\cite{DeligneWeil2},~(1.7.5); 
%  \cite{KatzSE},~(4.7.4)). It is an ascending filtration $V,$  indexed by $i\in \ZZ$
%    which satisfies  
%  $ N(M_i(V))\subset M_{i-2}(V),$ where  $N$ denotes  the logarithm of the inertial local monodromy (see \cite{DeligneWeil2},~(1.7.5)). This implies
  \begin{equation}\label{eqNFiltr}
  \Gr^M(\cJ_{n})=\bigoplus_{j=0}^{n-1} \bQl(-j),
  \end{equation}
where $M$ denotes the weight filtration (see \cite{DeligneWeil2},~(1.7.5)). Since multiplication with the logarithm of the usual topological generator of the tame inertia lowers the weight by $2$ (loc.~cit.), 
the Frobenius operation on the $I^t$-invariants is trivial. 

The defining properties of
$\cJ_n$ do not determine it as a $G^t$-module up to isomorphism if $n>1.$ Hence it would be more accurate to speak about a class of representations satisfying these properties. Since for us only the indecomposability and the triviality of the fixed space under the Frobenius operation is important, we neglect this ambiguity notationally. 

The restriction of $\cJ_n$ to some $G_l^t$ is  denoted by the same symbol.
\end{enumerate}

\begin{rem}\label{rembrauer} \begin{enumerate}
 \item Since ${\rm Gal}(\bar{k}/k)$ (resp. ${\rm Gal}(\bar{k}/k_l)$) is abelian
 any unramified irreducible representation of $G^t$ (resp. $G_l^t$) has rank one.

\item Any representation of rank one of $G^t$ (resp. $G_l^t$) is of the form 
$ V_{\chi'}\otimes F$ with $F$ an unramified rank-one representation of $G^t$ 
(resp. $G_l^t$), see~\cite{Laumon},~(3.5.3.1).  

\item Any {irreducible} $\bQl$-module $V$ of rank $l$ of $G^t$ is of the above used form
 $  \Ind_{G_l^t}^{G^t}( V_{\chi'}\otimes F)$  (loc.~cit.).  
 
 \item Let 
 $G_l^t\subseteq G_{l'}^t\subseteq G^t$ be a chain of proper inclusions such that
 $ V_{\chi'}$ is defined over the intermediate field $\FF_{q^{l'}},$ meaning that 
 the order of the image of the tame inertia divides $q^{l'}-1$. Then 
 $\Ind_{G_l^t}^{G_{l'}^t}( V_{\chi'}\otimes F)$ decomposes into several 
 rank-one factors since, by the operation of the Frobenius on the tame inertia determined by Formula~\eqref{eqzl1},
  the image of $G_{l'}^t$ under $\Ind_{G_l^t}^{G_{l'}^t}( V_{\chi'}\otimes F)$ is abelian. It hence follows from the functoriality of 
 the induction that $\Ind_{G_l^t}^{G^t}( V_{\chi'}\otimes F)$ is not irreducible. 
  \end{enumerate}   \end{rem}
%
% ,\quad \textrm{and hence}\quad   \Gr^M(\cJ_{n}\otimes \Ind_{G_l}^G(\cL_{\chi_l}\otimes F_l))=\bigoplus_{j=0}^{n-1} \Ind_{G_l}^G(\cL_{\chi_l}\otimes F_l))(-j)\,.
%\begin{rem}\label{remtamereps2} Recall from \cite{DeligneSommesTrig},~(4.14),~(4.15.1), that for two characters $\chi_1,\chi_2$ of $k^\times\,(k=\FF_q)$ which are not both trivial,
%and with $\chi_0:=\chi_1^{-1}\cdot \chi_2^{-1},$ one has associated the {\it Jacobi sum}
%$$ j(\chi_1,\chi_2):=-\sum_{\stackrel{x_1,x_2\in k^\times}{ x_1+x_2=-1}}\chi_1(x_1)\cdot \chi_2(x_2)$$ such that 
%\begin{equation}\label{eqjacobi} q\cdot j(\chi_1,\chi_2)=g(\chi_1^{-1}\chi_2^{-1},\psi)\cdot g(\chi_1,\psi)\cdot g(\chi_2,\psi)\,\end{equation}
%(in loc.cit., the Jacobi sum $j(\chi_1,\chi_2)$ is denoted $j((\chi_1\cdot \chi_2)^{-1},\chi_1,\chi_2)$ and the occurring Jacobi and Gau{\ss} sums are formed using 
%the inverse characters).  
%Then the Frobenius element of $G$ acts on the unramified 
%the $G$-module of rank one
%$$ J(\chi_1,\chi_2):=G(\chi_{1}^{-1}\chi_2^{-1},\psi)\otimes G(\chi_{1},\psi)\otimes G(\chi_2,\psi)(1)$$
%by the Jacobi sum  $j(\chi_1,\chi_2).$ 
% \end{rem}
 
% In the following, the symbol $\1$ stands for the trivial $\bQl^\times$ character of some 
% fundamental group~$G_l^t.$ 
% 
 
\begin{thm}\label{thmThom} For $k=\FF_q$ and $\ell\neq {\rm Char}(k)$ 
let  $F$ be an irreducible smooth $\bQl$-sheaf on $\AA^1_k\setminus S\stackrel{j}{\hookrightarrow}\AA^1_k,$ tamely ramified at $S\cup \infty,$ 
 such that $K=j_*F[1]\in {\rm Fourier}(\AA^1,\bQl).$ 
 Let $\cL_\chi$ be a nontrivial Kummer sheaf and suppose that $F$ is geometrically not 
 isomorphic to a translate 
 of the Kummer sheaf $\cL_{\chi^{-1}}.$ 
For fixed $s\in S$ with $k(s)=k_l,$  let  
$$ F_{\overline{\eta}_s}/F_{\overline{\eta}_s}^{I_s}=\bigoplus_{i=1}^{r_s}\cJ_{n_i}\otimes {\rm Ind}_{G_{l_i}^t}^{G_l^t} 
  ( V_{\chi_{i}}\otimes F_i)$$ be the decomposition into indecomposable 
  $G_l^t$-modules as in \eqref{eqNFiltr2}, where we identify $G_s^t$ 
  with $G_l^t$ and where $l\leq l_i.$  
  
Then $\MC_\chi(K)$ is an object in ${\rm Fourier}(\AA^1,\bQl)$ of the form $j_*H[1]$ with $H$ irreducible and smooth on $\AA^1\setminus S$   such one has a decomposition of $ H_{\overline{\eta}_s}/H_{\overline{\eta}_s}^{I_s}$ into indecomposable $G_l^t$-modules
  $$ H_{\overline{\eta}_s}/H_{\overline{\eta}_s}^{I_s}=\bigoplus_{i=1}^{r_s} H_i,$$
where $H_i$ is  as follows:
\begin{enumerate}
 \item If $\chi_i \neq \chi^{-1},\1$ then 
 $$ H_i=\cJ_{n_i}\otimes   {\rm Ind}_{G_{l_i}^t}^{G_l^t} 
\left(    V_{\chi\chi_i}\otimes ({F_{ \chi\chi_i(-1)}}\otimes G((\chi\chi_i)^{-1},\psi)\otimes G(\chi,\psi)\otimes G(\chi_i,\psi)\otimes F_i)\right)(1) \,,$$ 
  where  $F_{\chi\chi_i(-1)}$ stands for the   geometrically constant rank-one $G^t_{l^i}$-module whose ${\rm Frob}_{k_{l^i}}$-trace is $\chi\chi_i(-1)$ (note the convention of Rem.~\ref{rempullback} for the restriction of representations/sheaves, so that we interpret in the above formula $\cJ_{n_i}$ as the $G_s^t$-module $(\cJ_{n_i})_{\bar{\eta}_s}$ 
  and  $\chi$ as the character $\chi\circ {\rm Nm}^{k_{l^i}}_{k_l}$), and where the $(1)$ of the right hand side of the formula denotes a Tate twist. 
\item If $\chi_i =\1,$ then  $G_{l_i}^t=G_l^t$ and 
$$ H_i= \cJ_{n_i} \otimes  V_\chi  \otimes 
 F_i\,.$$ 

\item If $\chi_i=\chi^{-1},$  then  $G_{l_i}^t=G_l^t$ and 
$$ H_i = \cJ_{n_i}\otimes  F_{\chi(-1)} \otimes  F_i(-1)\,,$$ 
  where  $F_{\chi(-1)} $  stands for the 
  geometrically constant rank-one $G_s^t$-module whose ${\rm Frob}_{k_l}$-trace is $\chi(-1),$  and where the $(-1)$ of the right hand side of the formula denotes a Tate twist.  
\end{enumerate}

\end{thm}

 \begin{rem} \begin{enumerate}
 \item
For $k$ an algebraically 
 closed field this is proven in \cite{Katz96},~Cor.~3.3.6. For its proof 
 we use similar arguments, 
 further refined  by the results in Thm.~\ref{thmlaumon}. 
 \item It follows from Rem.~\ref{rembrauer}~(iii) that any $F$ with finite monodromy 
 satisfies the conditions on $ F_{\overline{\eta}_s}/F_{\overline{\eta}_s}^{I_s}$ of the above theorem. 
 \item In the dissertation of Tenzler \cite{Tenzler}, a variant of Thm.~\ref{thmThom} is used 
 which provides an algorithm to successively compute Frobenius traces for the construction 
 of rigid local systems 
 under the Katz algorithm \cite{Katz96}. As an application,  Frobenius traces 
 of certain $G_2$-motives are computed in \cite{Tenzler}. 
 \end{enumerate}
% \item Note that for any nontrivial $\chi,$ the Frobenius trace of $G(\chi,\psi)$ is given by the (negative of the) Gau{\ss} sum 
% $g(\chi,\psi)=-\sum_{a\in k^\times}\chi(a)\psi_k(a)$
%  (Thm.~\ref{thmlaumon}). Under the assumption $\chi_i \neq \chi^{-1},1,$ one has the well known relation 
%  $$ J(\chi,\chi_i)=-\frac{g(\chi,\psi)g(\chi_i,\psi)}{g(\chi\chi_i,\psi)}=\frac{-1}{q}\cdot (\chi\chi_i)(-1)\cdot g(\chi^{-1}\chi_i^{-1},\psi)g(\chi,\psi)g(\chi_i,\psi)\,,$$ 
%  where $J(\chi,\chi_i):=\sum_{a\in \FF_q}\chi(a)\chi_i(1-a)$ and using  
%  \begin{equation}\label{chichiinv} g(\chi,\psi)g(\chi^{-1},\psi)=\chi(-1)\cdot q\, ,\end{equation} cf.~\cite{DeligneSommesTrig}. 
%  Therefore (in Thm.~\ref{thmThom}~(i)) the 
%  trace of $\Frob_{q^{l_i}}$ on $\chi\chi_i(-\1)\otimes G(\chi^{-1}\chi_i^{-1},\psi)\otimes G(\chi,\psi)\otimes G(\chi_i,\psi)$ is equal to 
%  $-q^{l_i} J(\chi,\chi_i).$ Since $q^{l_i}$ is eliminated by the Tate twist, this implies that the trace of $\Frob_{q^{l_i}}$ on
%  $$\chi\chi_i(-\1)\otimes \cL_{\chi}\otimes \cL_{\chi_{i}}\otimes G(\chi^{-1}\chi_i^{-1},\psi)\otimes G(\chi,\psi)\otimes G(\chi_i,\psi)\otimes F_i(1)$$
%  is given by $-J(\chi,\chi_i)\cdot \tr(\Frob_{q^{l_i}},F_i).$ 
%             \end{enumerate}
\end{rem}
 
\proof 
Write  $\cF(K)=j'_*F'[1]\in {\rm Fourier}(\AA',\bQl)$ (where $\AA'$ is as in 
Section~\ref{secintrofourier}) for $F'$ the restriction of  $\cF(K)$ to $j':\GG_m'\hookrightarrow \AA'$
and let $$\cF(j_*\cL_{\chi}[1])=j'_*(\cL_{\chi^{-1}}\otimes G(\chi,\psi))[1]=j'_*H'[1]\in {\rm Fourier}(\AA',\bQl).$$ By Prop.~\ref{proprelfourier}
$$ \cF(\MC_\chi(K))=j_*(F'\otimes H')[1]\quad \Leftrightarrow\quad \MC_\chi(K)= a^*\cF'(j_*(F'\otimes H')[1])(1) $$ by Fourier inversion. 
Note that by our assumption on $F,$ the 
sheaf $j_*(F'\otimes H')[1]$ is  a Fourier sheaf which implies that $\MC_\chi(K)$ is again in ${\rm Fourier}(\AA,\bQl).$ 
Therefore Cor.~\ref{corstationaryback}  and Thm.~\ref{thmlaumon} imply that 
\begin{equation} \label{eqthmloc}  H_{\overline{\eta}_s}/H_{\overline{\eta}_s}^{I_s}= a^*\cF^{(\infty',0)}\left(\cF^{(0,\infty')}( F_{\overline{\eta}_s}/F_{\overline{\eta}_s}^{I_s})\otimes 
 V_{\chi}\otimes G(\chi,\psi)\right)(1)\,,
\end{equation}
where we write $ V_\chi$ instead of  $(\cL_\chi)^{-1}_{\overline{\eta}_{\infty'}}$ (resp.~$G(\chi,\psi)$ instead of $ G(\chi,\psi)_{\overline{\eta}_s}$) so that  
\begin{equation} \label{eqthmloc1} F_{\overline{\eta}_s}/F_{\overline{\eta}_s}^{I_s}=\bigoplus_{i}\cJ_{n_i}\otimes {\rm Ind}_{G_{l_i}}^{G} 
  ( V_{\chi_{i}}\otimes F_i)\end{equation}
by assumption. 
  We consider the individual  contribution of each summand $ \cJ_{n_i}\otimes {\rm Ind}_{G_{l_i}}^{G} 
  ( V_{\chi_{i}}\otimes F_i)$ of \eqref{eqthmloc1} to~\eqref{eqthmloc}:
  
 Assume first that $\chi_i \neq \chi^{-1},\1.$ Let 
  $$ H_i:=a^*\cF^{(\infty',0)}\left(\cF^{(0,\infty')}(\cJ_{n_i}\otimes {\rm Ind}_{G_{l_i}}^{G} 
  ( V_{\chi_{i}}\otimes F_i))\otimes 
 V_{\chi}\otimes G(\chi,\psi)\right)(1).$$ 
It follows from  the definition of $\cJ_n$    and from the exactness of the local Fourier transform  that tensoring tame irreducible $G_l^t$-modules with $\cJ_n$ commutes with local Fourier transformations (use an induction on $n,$ cf.~\cite{Tenzler}~Prop. 4.22(vi)). Hence 
\begin{eqnarray}\label{eqerste} \cF^{(0,\infty')}(\cJ_{n_i}\otimes {\rm Ind}_{G_{l_i}}^{G} 
  ( V_{\chi_{i}}\otimes F_i))&=&\cJ_{n_i}\otimes \cF^{(0,\infty')}({\rm Ind}_{G_{l_i}'}^{G'} 
  ( V_{\chi_{i}}\otimes  F_i))\nonumber \\
  &=&\cJ_{n_i}\otimes {\rm Ind}_{G_{l_i}'}^{G'} 
  ( V_{\chi_{i}}\otimes G(\chi_i,\psi)\otimes F_i)\nonumber
  \end{eqnarray} is again indecomposable.  Hence $H_i$ is indecomposable since 
  tensor products with rank-one sheaves and $\cF^{(\infty',0)}$  preserve indecomposability. 
  The push-pull formula for the induction and restriction of representations implies  
\begin{eqnarray}\nonumber H_i&=&a^* \cF^{(\infty',0)}\left(\cJ_{n_i}\otimes {\rm Ind}_{G_{l_i}'}^{G'} 
  \left[ V_{\chi_{i}}\otimes G(\chi_i,\psi)\otimes F_i\right]\otimes 
 V_{\chi}\otimes G(\chi,\psi)\right)(1)\\
&=& a^* \cF^{(\infty',0)}\left(\cJ_{n_i}\otimes {\rm Ind}_{G_{l_i}'}^{G'} 
  \left[ V_{\chi\chi_{i}}\otimes G(\chi_i,\psi)\otimes F_i\otimes G(\chi,\psi)\right]\right)(1)\,,\nonumber 
\end{eqnarray}
 where in the square bracket after the induction sign
we have restricted $  V_{\chi}\otimes G(\chi,\psi)$ to $G_{l_i}.$
This gives \begin{eqnarray}\nonumber H_i&=&
%a^*\cF^{(\infty',0)}\left(\cJ_{n_i}\otimes {\rm Ind}_{G_{l_i}'}^{G'} 
%  (\cL_{\chi_{i}}\otimes G(\chi_i,\psi)\otimes F_i ) \otimes \cL_{\chi}\otimes G(\chi,\psi))\right)(1) \nonumber \\
%  &=&
   a^*\left(  \cJ_{n_i}\otimes  {\rm Ind}_{G_{l_i}'}^{G'} 
  \left[ V_{\chi\chi_{i}}\otimes G(\chi^{-1}\chi_i^{-1},\psi)\otimes G(\chi,\psi)\otimes G(\chi_i,\psi)\otimes F_i)\right]\right)(1)\\
  &=&
 \cJ_{n_i}\otimes  {\rm Ind}_{G_{l_i}'}^{G'} 
  ( V_{\chi\chi_{i}}\otimes F_{ \chi\chi_i(-1)}\otimes G(\chi^{-1}\chi_i^{-1},\psi)\otimes G(\chi,\psi)\otimes G(\chi_i,\psi)\otimes F_i)(1)\,,
  \nonumber 
  \end{eqnarray} using Thm.~\ref{thmlaumon}(iii).
   Note 
  for the last equation, via the trace function of Kummer sheaves (cf.~Section~\ref{subsub2}), 
 the effect of $a^*$ on the associated Frobenius trace in the above formula for $H_i$ amounts to a multiplication 
 with $\chi\chi_i(-1).$

 If $\chi_i=\1,$ then we have $G_{l_i}^t=G_l^t$ by indecomposability
 (see Rem.~\ref{rembrauer} (iv)) and Thm.~\ref{thmlaumon} implies that 
 $$ \cF^{(0,\infty')}(\cJ_{n_i}\otimes  F_i)=\cJ_{n_i}\otimes  F_i.$$ Hence, 
 with $$H_i:=  a^* \cF^{(\infty',0)}\left(\cF^{(0,\infty')}( \cJ_{n_i}\otimes F_i)\otimes 
 V_{\chi}\otimes G(\chi,\psi)\right)(1)$$ we obtain from Thm.~\ref{thmlaumon}
\begin{eqnarray}\nonumber H_i&=&
a^*\cF^{(\infty',0)}\left(\cJ_{n_i}\otimes  V_{\chi}\otimes G(\chi,\psi)\otimes F_i)\right)(1)\nonumber \\
  &=& a^* (\cJ_{n_i}\otimes  V_\chi\otimes  G(\chi^{-1},\psi)\otimes G(\chi,\psi)\otimes F_i(1))\nonumber \\
  &=&\cJ_{n_i} \otimes  V_\chi \otimes F_i\,,
  \end{eqnarray}
  where we used $\chi=\chi_i^{-1}$ and Formula~\eqref{chichiinv}.   
  
  If $\chi_i=\chi^{-1},$ then again have  $G_{l_i}^t=G_l^t$ and 
 with $$H_i:=  a^* \cF^{(\infty',0)}\left(\cF^{(0,\infty')}( \cJ_{n_i}\otimes  V_{\chi^{-1}}\otimes F_i)\otimes 
 V_{\chi}\otimes G(\chi,\psi)\right)(1)$$ we obtain 
\begin{eqnarray}\nonumber H_i&=&
a^*\cF^{(\infty',0)}\left(\cJ_{n_i}\otimes   V_{\chi^{-1}}\otimes G(\chi^{-1},\psi)\otimes V_{\chi}\otimes G(\chi,\psi)\otimes F_i)\right)(1) \nonumber \\
  &=&  \cJ_{n_i}\otimes F_{\chi(-1)}\otimes  F_i(-1)\,,\nonumber 
  \end{eqnarray}
  where we used the identity \eqref{chichiinv} as well as Thm.~\ref{thmlaumon}, (ii), (iv) and~(v).
\Endproof

\section{The determinant of the \'etale middle convolution.}\label{sec:det}

\subsection{Local epsilon constants, local Fourier transform,  and Frobenius determinants.}\label{sec4.1}
As in the previous sections, we fix a finite field $k=\FF_q\, (q=p^m)$ and an  additive $\overline{\mathbb{\QQ}}_\ell^\times$-character 
$\psi$ of $\AA^1(\FF_q)\, (\ell\neq p).$  By composing with the trace function, the character $\psi$ uniquely determines
additive $\bQl^\times$-characters of any extension 
$\AA^1(\FF_{q^l}),$ again denoted $\psi.$ \\

Recall the theory of local epsilon constants for $X$ a connected smooth projective curve over~$k,$  see 
\cite{DeligneFonctionelles} and \cite{Laumon}, Section~3:
 The 
 $L$-function of $K\in \rD^b_c(X,\bQl)$ is defined as 
 $$ L(X,K;t)=\prod_{x\in |X|}\frac{1}{\det(1-t^{\deg x}\cdot \Frob_x,K)}.$$ By the work of Grothendieck,
this $L$-function is the product expansion of 
 $$ \det(1-t.\Frob_q,R\Gamma(X\otimes_{{k}} \overline{k},K))^{-1}$$ and it satisfies the following functional 
 equation:
 \begin{equation}\label{eql2}
 L(X,K;t)=\epsilon(X,K)\cdot t^{a(X,K)}\cdot L(X,D(K)),\end{equation}
 where $D(K)$ denotes the Verdier dual of $K$ and where $a(X,K)$ and 
 $\epsilon(X,K)$ are defined as follows: 
 $$ a(X,K):=-\chi(X,K) \quad \textrm{(Euler characteristic as defined in~ loc.cit., Section~0.8)}$$
 \begin{equation}\label{eql12} 
 \epsilon(X,K):=\det(-\Frob_q,R\Gamma(X_{\overline{k}},K))^{-1},\end{equation}
  %Then one has for each distinguished triangle $K'\to K\to K''\to K'[1]$ in $\Perv(X,\bQl)$ the following equalities
  
By the work of Deligne \cite{DeligneFonctionelles}, there is a unique map $\epsilon$ 
which, depending on a fixed 
character $\psi$ as above, associates to a triple $(T,K,\omega)$  
($T=(T,\eta,t)$ a henselian trait with residue field $k,$ $K\in \rD^b_c(T,\bQl),$  $\omega$ a nontrivial meromorphic $1$-form on $T$) a
{\it local epsilon constant} $\epsilon(T,K,\omega)\in \bQl^\times$ such that the following axioms hold (see \cite{DeligneFonctionelles}, Thm. 4.1,  and the reformulation in \cite{Laumon},~Thm.~(3.1.5.4)): 

 \begin{prop} \label{remdeligneaxioms} {\it 
  \begin{enumerate}
 \item The association $(T,K, \omega )\mapsto \epsilon(T,K, \omega )$ depends only on 
 the isomorphism class of the triple $(T,K, \omega ).$
 \item For any distinguished triangle $K'\to K \to K'' \to K'[1]$ in $\rD^b_c(T,\bQl)$ one has 
 \begin{equation}\label{eql44} 
 \epsilon(T,K, \omega )=\epsilon(T,K', \omega )\cdot\epsilon(T,K'', \omega ).
 \end{equation}
 \item If $K$ is supported on the closed point $t$ of $T$ 
 then \begin{equation}\label{eqlaumona} \epsilon(T,K,\omega)=\det(-\Frob_t,K)^{-1}.\end{equation}
 \item If $\eta$ denotes the generic point of $T,$ if $\eta_1/\eta$ is a finite separable extension of $\eta$ and if $f:T_1\to T$ denotes the 
 normalization of $T$ inside $\eta_1,$ then for any $K_1\in \rD^b_c(T_1,\bQl)$ such that 
 the generic rank $r(K_1)$ of $K_1$ (extended from the usual generic rank by additivity, so that in odd cohomological degree the rank is counted negative) is equal to $0$ (\cite{DeligneFonctionelles}), one has 
 \begin{equation} \label{eql5}
 \epsilon(T,f_*K_1,\omega)=\epsilon(T_1,K_1,f^*\omega).
 \end{equation} (The condition on the rank is rather mild because one can replace  $K_1$ 
 by $\tilde{K}_1:=K_1\oplus (\bQl^{e_1}[e_2])$ such that $r(\tilde{K}_1)=0,$ see \cite{DeligneFonctionelles} and \cite{Laumon},~(3.3.2).)
 
 \item If $V$ denotes a rank-one local system on $\eta$ corresponding to a character $\mu:
 K_x^\times \to \overline{\mathbb{\QQ}}_\ell^\times$ via reciprocity and if $j:\eta \hookrightarrow T$ denotes 
 the obvious inclusion, 
 then $ \epsilon(T,j_*V,\omega) $ coincides with Tate's local constant associated to $\mu$ (cf.~\cite{Laumon},~(3.1.3.2)). \end{enumerate}}
 \end{prop}

 For $V$ a smooth sheaf on $\eta$ and $j:\eta\hookrightarrow T$ the inclusion  
 (with $V$ viewed as a $G$-module via the fibre functor at $\bar{\eta}$ and with inertial invariants $V^I\simeq (j_*V)_{\bar{t}}$) one defines 
  $$\epsilon(T,V,\omega):=\epsilon(T,j_*V,\omega) $$
 and 
 \begin{equation}\label{eqepsilonnull}\epsilon_0(T,V,\omega):=\epsilon(T,j_!V,\omega)=
 \epsilon(T,V,\omega)\cdot \det(-\Frob_t,V^I),\end{equation}
 where the last equality follows from Prop.~\ref{remdeligneaxioms}(ii),(iii), using that in 
 $D_b^c(T),$ the extension by zero $j_!V$ is an extension of $j_*V$ (placed in cohomological degree zero) and the shifted invariants $V^I[-1] ,$ supported at the closed point $t,$
 cf.~\cite{DeligneFonctionelles}, Chapitre~(5).\\

 If $x$ is a closed point of $X,$ then $X_{(x)}$ denotes the Henselization of $X$ at $x$ (cf.~\cite{Laumon},~Section~0.4).  
 By the work of Laumon (\cite{Laumon}, Thm. 3.2.1.1), the epsilon constant decomposes into a product of local epsilon constants, depending on a  
 nontrivial 
 meromorphic differential
 $1$-form $\omega$ on $X,$
 as follows:
 \begin{equation}\label{eql3}
 \epsilon(X,K)=q^{C(1-g(X))r(K)}\prod_{x\in |X|}\epsilon(X_{(x)},K|_{X_{(x)}}, \omega |_{X_{(x)}}),\end{equation}
 where $C$ denotes the cardinality of connected components of $X\times \bar{k}$ and where $g(X)$ is the genus of some component of 
 $X\times \bar{k}.$ Another variant of \eqref{eql3} is 
 \begin{equation}\label{eql4}
 \epsilon(X,j_*F)=q^{C(1-g(X))r(F)}\prod_{x\in |X|}\epsilon(X_{(x)},F|_{\eta_{x}, }, \omega |_{X_{(x)}}),\end{equation}
 where $F$ is smooth on on an dense open $U\subseteq X$ with $j:U\hookrightarrow X$
 (\cite{Laumon},  (3.2.1.6)).\\

 The  local epsilon constants satisfy the following additional properties:

 \begin{prop}\label{reml1}(Laumon){\it  \begin{enumerate} \item 
 If 
 $$ 0\to V'\to V\to V''\to 0$$ is a short exact sequence of $G$-modules then 
\begin{equation}\label{eqewllshortex}
 \epsilon_0(T,V,\omega)=\epsilon_0(T,V',\omega)\cdot \epsilon_0(T,V'',\omega),
\end{equation}
see \cite{Laumon},~(3.1.5.7).   
%  As a special case of Eq.~\eqref{eql3},  let $U\stackrel{j}{\hookrightarrow} \PP^1_k$  be a dense  open subscheme
%with $S=\PP^1\setminus U$ and let $F$ be a smooth sheaf on $U$ of rank $r>0$ which is smooth at $\infty.$ 
%For $s\in |X|$  and  $\omega_0=-dx$ (with $x$ denoting the affine coordinate of $\AA^1\subset \PP^1$) so that  
%$$ \epsilon_0(X_{(s)},F_{\eta_s},\omega_0|_{X_{(s)}}):=\epsilon(X_{(s)},j_!(F|_{\eta_s}),\omega_0|_{X_{(s)}})$$ 
%with $j$ denoting the inclusion of $\eta_s$ into $X_{(s)}.$ 
%Then \cite{Laumon},~Thm.~3.3.1.2, states that 
% \begin{equation}\label{eqLaumonepsilonproduit}
%  {\rm det}(-\Frob_q,R\Gamma_c(U\otimes_k \overline{k},F))^{-1}\cdot q^r\cdot {\rm det}(-\Frob_\infty,F_{\overline{\infty}})=\prod_{s\in S}\epsilon_0(X_{(s)},
% F|_{{\eta}_s}, {\omega_0}|_{X_{(s)}}) \,.
% \end{equation}
 \item 
%                           \item Let $ V_\chi$ be the Kummer sheaf on $T$ 
%                           attached to $\chi:k^\times \to \bQl$ (viewed as a $G$-module) and let
%                           $F$ be an unramified $G$-module of rank one.                      Then
%                       $$ \epsilon_0(T, V_\chi\otimes F,d\pi)=-\det(\Frob_t,F)\cdot 
%                       \epsilon_0(T,L_\chi,d\pi),$$ see \cite{Laumon} (3.1.5.6), (3.5.3.1).
                           Let $K_t$ denote the completion of the function field of the generic point $\eta$ of 
                           $T$ and let $\nu_t:K_t^\times \to \ZZ$ its natural valuation. 
                           By \cite{Laumon}, 3.1.5.6, if $K\in \rD^b_c(T,\bQl)$ and if $F$ is a smooth sheaf on $T$ then 
 \begin{equation}\label{eql7}
 \epsilon( T, K\otimes F,\omega)=\epsilon(T,K,\omega)^{r(F)}\cdot \det(\Frob_t,F)^{a(T,K,\omega)},
 \end{equation} where $a(T,K,\omega)$ is defined as follows (cf.~\cite{Laumon},~(3.1.5.1), (3.1.5.2)):
 $$ a(T,K,\omega)=r(K_{\bar{\eta}})+s(K_{\bar{\eta}})-r(K_{\bar{t}})+r(K_{\bar{\eta}})\nu_t(\omega),$$
 where $s(K_{\bar{\eta}})$ is the Swan conductor of $K_{\bar{\eta}}$ 
 (which vanishes if and only if $K_{\bar{\eta}}$ is tame, cf.~\cite{Laumon}, (2.1.4))
and where $\nu_t(a\cdot db)=\nu_t(a)$ for $a\cdot db\in \Omega^1_{K_t}\setminus 0$ and $\nu_t(b)=1.$   
 \item Let $k_1$ be a finite extension of $k.$ Let $V$ be an 
   irreducible $G$-module of the form $f_* V_1$ with $f:T_1=T\otimes_k k_1\to T$ and 
with $ V_1$ tame. Let $G_1=\Gal(\bar{\eta}/\eta_1)$ where $\eta_1$ denotes the generic 
point of $T_1.$ Then
\begin{equation}\label{eqepsind}
 \epsilon_0(T,V,d\pi)=\epsilon_0(T_1, V_1,d\pi_1),
\end{equation}
where $\pi_1$ is a uniformizer of $T_1$ induced by $\pi,$ cf.~\cite{Laumon},~3.5.3.1.      
\item 
If the character $\chi$ is nontrivial, then 
  $ j_! V_\chi=j_* V_\chi.$ Hence  if $K=j_* V_\chi$ then 
  \begin{equation}\label{eql8} \epsilon(T,K,d\pi)=\epsilon_0(T, V_\chi, d\pi)=\chi(-1) g(\chi,\psi)
  \end{equation}
with $g(\chi,\psi)$ the
  Gau{\ss} sum occurring in  Prop.~\ref{propkatz1}
  (\cite{Laumon},~(3.5.3.1)).   If $\chi$ is trivial then  also 
  $$\epsilon_0(T, V_\chi, d\pi)=-1\,\quad \textit{and}\quad \epsilon(T, V_\chi, d\pi)=1$$
%  Moreover,
%  $$ \epsilon(T, V_\chi,-d\pi)=-\epsilon(T, V_\chi,d\pi)\textrm{ and }
%  \epsilon_0(T, V_\chi,-d\pi)=-\epsilon_0(T, V_\chi,d\pi),$$
%  see \cite{DeligneFonctionelles}, Thme.~4.1. 
  \item For $a\in k(\eta)^\times $ one has 
  $$ \epsilon(T,K,a\omega)=\chi_K(a)q^{r(K_{\overline{\eta}})\nu_t(a)}\epsilon(T,K,\omega)\,$$ 
  where $\chi_K:K_t^\times\to \bQl^\times$ is the character induced by  $\det(K_\eta)$ via reciprocity (\cite{Laumon}, (3.1.5.5)).
  \item The behaviour of local epsilon constants under Tate twists is given as follows
  (\cite{Laumon}~(3.2.1.4):
  $$ \epsilon(X_{(x)},K(m)|_{X_{(x)}},\omega|_{X_{(x)}})=q_x^{-m\cdot a(X_{(x)},K|_{X_{(x)}},\omega|_{X_{(x)}})}\epsilon(X_{(x)},K|_{X_{(x)}},\omega|_{X_{(x)}}),$$
%  where $K(m)$ denotes the $m$-th Tate twist of $K$ and where $q_x=q^{\deg(x)}.$ Especially, for a smooth rank-one sheaf $F$ on $T,$ 
%  $$ \epsilon(T,F(m),\omega)=q^{-m}\epsilon(T,F,\omega)\,.$$
 % \item If $k=\FF_q$ and $X,K,\omega,S,x$ are as above and if $X_m,K_m,\omega_m,S_m,x_m$ are the objects
 % which are obtained from a respective basechange
 % to $\FF_{q^m}\,(m\in \NN_{>0}),$ then it is proved in \cite{Laumon}, Preuve de (3.3.2.1), that 
 % $$\prod_{x_m\in |X_m|}\epsilon({X_m}_{(x_m)},K_m|{X_m}_{(x_m)},
 % {\omega_m} |_{{X_m}_{(x_m)}})=(-1)^{(m-1)a(X,K)} \prod_{x\in |X|}\epsilon(X_{(x)},K|X_{(x)}, \omega |_{X_{(x)}})^m\,,$$
%where $ a(X,K)=-\chi(X,K).$ 
%
%Especially, with  the notation and assumptions of (i) above with $K=j_!F,$  $F_m=F|_{U_m}$ and with 
%$\eta_{m,s}$ denoting the generic point of ${X_m}_{(s)},$ then 
%$$\prod_{s_m\in S_m}\epsilon_0(X_{m}_{(s_m)},F_m|_{\eta_{s_m}}, {\omega_{0,m}} |_{X_{m}_{(s_m)}})=
%(-1)^{(m-1)(a(X,K)+r)} \prod_{s\in S}\epsilon_0(X_{(s)},F_{\eta_s}, \omega_0 |_{X_{(s)}})^m\,.$$
\end{enumerate}}\Endproof
 \end{prop}
 
% \begin{rem}\label{remextension} Let $T=(T,\eta,t)$ be a henselian trait as above with 
% residue field $k=\FF_q$ and Galois group $G=\pi_1(\eta,\bar{\eta}).$  Let  $V$ be an irreducible 
% $G^t$-module of rank $l$ of $G^t$ of the  form
% $  V=\Ind_{G_l^t}^{G^t}( V_{\chi'}\otimes F)$  as in Rem.~\ref{rembrauer}, where $G_l=\pi_1(\eta_l,\bar{\eta}_l)$ is the Galois group of the trait $$T_l=T\otimes_k k_l=
% (T_l,\eta_l,t_l)\,\, (k_l=\FF_{q^l}).$$  Then $f:T_l\to T$ is unramified of rank $l.$ If $\pi_l$ the uniformizer of $ T_l$ induced by $\pi,$ then 
% $$\epsilon_0(T,V,d\pi)=\epsilon_0(T_l, V_{\chi'}\otimes F,d\pi_l)$$ (see \cite{Laumon}~(3.5.3.1) for a proof, using Prop.~\ref{remdeligneaxioms}(iv)). It follows hence from 
% 
% \end{rem}

\begin{prop}  \label{propindecomp} Let $G=\pi_1(\eta,\bar{\eta})$ be the Galois group of the generic point of a henselian trait 
$T=(T,\eta,t)\,, t=\Spec(k),k=\FF_q,$ which uniformizer $\pi$  as above. Let 
$$T_l=T\otimes_k k_l=(T_l,\eta_l,t_l)\, (k_l=\FF_{q^l})$$ and $G_l=\pi_1(\eta_l,\bar{\eta}_l).$ 
Suppose that  $V$ is a tame indecomposable $G$-module  which, in the notation of  Eq.~\eqref{eqNFiltr2}, can be written as 
   \begin{equation}\label{eqNFiltr22} V=
   \cJ_{n}\otimes \Ind_{G_l^t}^{G^t}( V_{\chi_l}\otimes F_l)\,,
   \end{equation}
   where $ V_{\chi_l}$ is the $G_l$-module derived from the Kummer sheaf attached to $\chi_l:\GG_m(k_l)\to \bQl^\times$ (see Thm.~\ref{thmlaumon}(iii))
   and where $F_l$ is an unramified $G_l^t$-module. 
Then the following holds:
\begin{enumerate}
\item  If $\chi\neq \1$ then $$ \epsilon(T,V,d\pi)=\epsilon_0(T,V,d\pi)=  
q^{n(n-1)/2}\cdot \left(\chi_l(-1) g(\chi_l,\psi)\cdot \det(\Frob_{t_l}, F_l)\right)^n \cdot  \,.$$
\item If $\chi= \1$ then $l=1$ and we can write 
$$ V=
   \cJ_{n}\otimes F$$ with $F$ an unramified $G^t$-module.  Then 
$$ \epsilon_0(T,V,\omega)=(-q)^{n(n-1)/2}\det(-\Frob_t,F)^n$$ and 
  $$ \epsilon(T,V,\omega)=(-q)^{n(n-1)/2}\det(-\Frob_t,F)^{(n-1)}=\epsilon_0(T,V/V^I,\omega).$$
\end{enumerate}
\end{prop}

\proof Let us first treat the case where ${\chi_l}\neq \1:$ By definition of $\cJ_n$ (see Eq.~\eqref{eqNFiltr}),
$$
  \Gr^M(V)=\bigoplus_{j=0}^{n-1} \Ind_{G_l}^G( V_{\chi_l}\otimes F_l)(-j)\,.$$
  Therefore Prop.~\ref{reml1}~(i) implies that 
  \begin{eqnarray}\label{eqlaumonb}
  \epsilon_0(T,V,d\pi)&=&\prod_{j=0}^{n-1}\epsilon_0(T,\Ind_{G_l}^G( V_{\chi_l}
  \otimes F_l)(-j),d\pi)\,.
  \end{eqnarray}
  Since ${\chi_l}$ is nontrivial,
$$j_*(\Ind_{G_l}^G(  V_{\chi_l}\otimes F_l))=j_!(\Ind_{G_l}^G( V_{\chi_l}\otimes F_l)),$$ and
therefore
$$\epsilon(T,\Ind_{G_l}^G(  V_{\chi_l}\otimes F_l),d\pi)=\epsilon_0(T,\Ind_{G_l}^G(  V_{\chi_l}\otimes F_l),d\pi),$$
as well as 
$$\epsilon_0(T,V,d\pi)=\epsilon(T,V,d\pi).$$
%where $\pi_l$ denotes the uniformizer of $R_l$ (where $R_l$ is the henselian ring underlying 
%$T_l=T\times_k k_l$).
It follows from Prop.~\ref{reml1}(iv) that 
\begin{equation}\label{eqlaumonb2}
\epsilon_0(T,\Ind_{G_l}^G(  V_{\chi_l}\otimes F_l),d\pi)=
\epsilon_0(T_l,  V_{\chi_l}\otimes F_l,d\pi_l).
\end{equation} Using Prop.~\ref{reml1}(ii), especially Eq.~\eqref{eql7}, for the first equality
(using that the valuation of $d\pi_l$ is zero and $F_l$ is smooth on $T_l$), and Prop.~\ref{reml1}(iv) for the second equality, one obtains 
\begin{equation}\label{eqlaumonb22}
\epsilon_0(T_l,  V_{\chi_l}\otimes F_l,d\pi_l)=\epsilon_0(T_l,  V_{{\chi_l}},d\pi_l)\cdot \det(\Frob_{t_l}, F_l)=\chi_l(-1) g(\chi_l,\psi)\cdot \det(\Frob_{t_l}, F_l)\,.
\end{equation}

Therefore, taking Tate twists into account, 
   \begin{eqnarray} 
  \epsilon_0(T,V,d\pi)=  \epsilon(T,V,d\pi)
  &=& \prod_{j=0}^{n-1}q^{j}\epsilon(T,\Ind_{G_l}^G( V_{\chi_l}\otimes F_l(-j)),d\pi)\label{labeleq0}\\
  &=&q^{n(n-1)/2}\left(\chi_l(-1) g(\chi_l,\psi)\cdot \det(\Frob_{t_l}, F_l)\right)^n\,.
  \end{eqnarray}

  Let now $\chi=\1:$ it follows from indecomposability that $G_l=G,$ see
  Rem.~\ref{rembrauer}(iv). Hence we have 
  $$ V=  \cJ_{n}\otimes  F$$ 
  with $F$ smooth on $T$ of rank one. As above we obtain
   \begin{eqnarray}\label{eqlaumonb3}
  \epsilon_0(T,V,d\pi)&=&\prod_{j=0}^{n-1}\epsilon_0(T,F(-j),d\pi)\,
  \end{eqnarray} since $\epsilon_0$ is multiplicative on short exact sequences of $G$-modules.
It follows from Prop.~\ref{reml1}(ii),(iv) that $$\epsilon(T,F(-j),d\pi)=\epsilon(T,\bQl\otimes F(-j),d\pi)=
\epsilon(T, V_{\1},d\pi)^1\cdot \det(\Frob_t,F(-j))^0=1$$
since the valuation of $d\pi$ is zero and $F(-j)$ is smooth (note that $ V_{\1}=\bQl$). 
Since, by Eq.~\eqref{eqepsilonnull},  for general $G$-modules  $W$
  $$\epsilon_0(T,W,\omega)=
 \epsilon(T,W,\omega)\cdot \det(-\Frob_t,W^I),$$ we obtain 
 $$\epsilon_0(T,F(-j),d\pi)= \det(-\Frob_t,F(-j))=(-q)^j\det(-\Frob_t,F).$$
 Consequently,
$$ \epsilon_0(T,V,\omega)=(-q)^{n(n-1)/2}\det(-\Frob_t,F)^n$$ and 
  $$ \epsilon(T,V,\omega)=(-q)^{n(n-1)/2}\det(-\Frob_t,F)^{(n-1)}.$$ The last equation follows 
  from $$\epsilon_0(T,V,\omega)=\epsilon_0(T,V^I,\omega)\epsilon_0(T,V/V^I,\omega)$$
 (see Prop.~\ref{reml1}(i)) and $\epsilon_0(T,V^I,\omega)=\det(-\Frob_t,F).$
  \Endproof

 \subsection{The Frobenius determinant of the middle convolution.}
 Let $k=\FF_q\,(q=p^a),$ let $U\stackrel{j}{\hookrightarrow} \AA^1_k$  be a dense  open subscheme, and let 
$S=\AA^1_k\setminus U.$  Let $V=j_*F[1]\in {\rm Fourier}(\AA^1_k,\bQl)$  be an 
irreducible middle extension sheaf, where $F$ is an irreducible  smooth sheaf of rank $r(F)$ on $U$ which is
tamely ramified in $S\cup \infty.$
Consider the following conditions on   $F$   with respect to a nontrivial character $\chi:\GG_m(k)\to \bQl^\times:$ 

\begin{enumerate}
%\item $F$ is an irreducible  nontrivial smooth sheaf of rank $r(F)$ on $U$ such that   $V:=j_*F[1]$ is an
%object in ${\rm Fourier}(\AA^1_k,\bQl).$
%\item   $F$  is
%tamely ramified in $S\cup \infty.$
\item  $F$ has 
scalar inertial  local  monodromy
at $\infty,$ whose restriction to the tame inertia group is given by the (restriction of the)  Kummer sheaf associated to a character $\chi:k^\times \to \bQl,$
see~Section~\ref{subsub2}.
\item   $F$ is geometrically not isomorphic to a translated Kummer sheaf of the form $\cL_{\chi^{-1}}(y-x),$
where $x$ denotes the coordinate of $\AA^1$ and $y\in | \AA^1|.$ 
\end{enumerate}

\begin{defn}{\rm  If $F$ as above satisfies the above conditions (i),(ii) we say that  $F$ is in {\it standard form} with respect  to $\chi.$ }
\end{defn}

\begin{rem}\label{remstandardform} We remark that if  $F,$ is in standard form with respect  to $\chi,$ then 
$\MC_\chi(V)$ is in standard form with respect  to $\chi^{-1}$ (see \cite{Katz96},~Cor.~3.3.6). \end{rem}

 For $y\in \AA^1(k)\setminus S,$ let $U_y=U \setminus y,$ and let 
  $${F_y(\chi)}=F|_{U_y}\otimes \cL_{\chi}(y-x)|_{U_y}.$$    
Note that $j_{*}{F_y(\chi)}$ is smooth in $\infty,$ where $j:U_y\hookrightarrow \PP^1$ is the natural inclusion. We define ${F_y(\chi)}_{\bar{\infty}}:=j_*{F_y(\chi)}_{\bar{\infty}}.$ 
We assume further that 
for each $s\in S$ with $k(s)=k_{l_s}$ one has a decomposition of   
$ F_{\overline{\eta}_s}/F_{\overline{\eta}_s}^{I_s^t}$ into indecomposable $G_s^t=G_{l_s}^t$ modules (where $G_s^t$ is the tame quotient of 
$G_s=\pi_1(\overline{\eta_s}/\eta_s)$ with tame inertia subgroup $I_s^t$ and where
$G_{l_s}=\pi_1(\overline{\eta}_{l_s}/\eta_{l_s})$ with $\eta_{l_s}$ the generic point of 
$T_{l_s}=T\otimes_k k_{l_s}$ with $T=\Spec(k\{\pi\})$ the trait of the henselization of $k[\pi]_{(\pi)}$)
\begin{equation}\label{eqassumptiononv1} F_{\overline{\eta}_s}/F_{\overline{\eta}_s}^{I_s^t}\simeq \bigoplus_{i_s=1}^{r_s}  \cJ_{n_{i_s}}\otimes \Ind_{G_{l_{i_s}}^t}^{G_{l_s}^t}( V_{\chi_{i_s}}\otimes F_{i_s})\,,\end{equation}
 with $F_{i_s}$ unramified (using the notation of Section~\ref{secThom}). Since the sheaf $\cL_{\chi}(y-x)$
 is smooth at $s,$ the  latter formula implies that  ${F_y(\chi)}_{\overline{\eta}_s}/{F_y(\chi)}_{\overline{\eta}_s}^{I_s^t}$ is of the form
    %Sometimes we write $\Frob_{k(x)}$ instead of $\Frob_q^{\deg(x)}$ and $q_x=q^{\deg(x)}.$  
\begin{equation}\label{eqassumptiononv} {F_y(\chi)}_{\overline{\eta}_s}/{F_y(\chi)}_{\overline{\eta}_s}^{I_s^t}\simeq \bigoplus_{i_s=1}^{r_s}  \cJ_{n_{i_s}}\otimes \Ind_{G_{l_{i_s}}^t}^{G_{l_s}^t}( V_{\chi_{i_s}}\otimes F_{i_s}')\,\end{equation}
with $F_{i_s}'$ unramified.

\begin{thm}\label{thml2} Under the above  assumptions, let $F$ be in standard form with respect $\chi$ and 
  %Sometimes we write $\Frob_{k(x)}$ instead of $\Frob_q^{\deg(x)}$ and $q_x=q^{\deg(x)}.$  
%Suppose that   $F$  is a smooth sheaf of rank $r(F)$ on $U$ satisfying the conditions above  which is in standard form with respect to $\chi.$  
    let  $\omega_0:=dx.$ Then
    % let $j:U_y\hookrightarrow \PP^1\,(U_y=U\setminus y)$ denote the natural inclusion,
%     and let ${F_y(\chi)}=F|_{U_y}\otimes \cL_{\chi}(y-x)|_{U_y}.$  Then the the following holds:  
$$ \det(-\Frob_q,H^1(\PP^1_{\overline{k}},j_*{F_y(\chi)})) =q^{-r(F)}
 \det(\Frob_{\infty}, {F_y(\chi)}_{\bar{\infty}})^{-2} 
\prod_{s\in S\cup y}\epsilon_0(\PP^1_{(s)},
 {F_y(\chi)}_{{\eta}_s}/{F_y(\chi)}_{{\eta}_s}^{I_s^t}, \omega_0|_{\PP^1_{(s)}})\,,$$ 
 where for $s=y,$ 
 $$  \epsilon_0(\PP^1_{(y)},
 {F_y(\chi)}_{{\eta}_y}/{F_y(\chi)}_{{\eta}_y}^{I_y}, \omega_0|_{\PP^1_{(y)}})=
 \epsilon_0(\PP^1_{(y)},
 {F_y(\chi)}_{{\eta}_y}, \omega_0|_{\PP^1_{(y)}})=
 (\chi(-1)g(\chi,\psi))^{r(F)}\cdot \det(\Frob_y,{F_y(\chi)}_{\bar{\eta}_y}),$$
 and where for  $s\in S,$ 
 $$  \epsilon_0(\PP^1_{({s})},
 {F_y(\chi)}_{{\eta}_s}/{F_y(\chi)}_{{\eta}_s}^{I_s^t}, \omega_0|_{\PP^1_{({s})}})=
 \prod_{{i_{s}=1}}^{r_s}\epsilon_0( \PP^1_{({s})},
 \cJ_{n_{i_s}}\otimes \Ind_{G_{s,l_{i_s}}}^{G_s}( V_{\chi_{i_s}}\otimes {F_{i_s}'}), \omega_0|_{\PP^1_{({s})}})\,,$$
 with 
 $$\epsilon_0( \PP^1_{({s})},
 \cJ_{n_{i_s}}\otimes \Ind_{G_{s,l_{i_s}}}^{G_s}( V_{\chi_{i_s}}\otimes {F_{i_s}'}), \omega_0|_{\PP^1_{({s})}})=\quad \quad \quad \quad \quad \quad\quad \quad \quad\quad \quad \quad \quad \quad \quad\quad \quad \quad \quad \quad\quad \quad \quad$$
 $$\quad \quad \quad\quad \quad \quad \quad \quad 
  \left\{
 \begin{array}{ll}
q^{l_{i_s}\cdot n_{i_s}(n_{i_s}-1)/2)}
\cdot (\chi_{l_{i_s}}(-1) g(\chi_{l_{i_s}},\psi)\cdot
 \det(\Frob_{t_{l_{i_s}}}, F_{l_{i_s}}'))^{n_{i_s}}& \quad \textrm{if}\quad  \chi_{l_{i_s}} \neq \1 \\
 (-q)^{l\cdot n_{i_s}(n_{i_s}-1)/2}\cdot \det(-\Frob_s,F_{i_s}'))^{n_{i_s}}& \quad \textrm{if}\quad  \chi_{l_{i_s}} =\1\,.
  \end{array}\right. $$
% 
% $$\epsilon_0( \PP^1_{({s})},
% \cJ_{n_{i_s}}\otimes \Ind_{G_{s,l_{i_s}}}^{G_s}( V_{\chi_{i_s}}\otimes {F_{i_s}'}), \omega_0|_{\PP^1_{({s})}})=\begin{array{c}
%q^{l_{i_s}}\cdot (n_{i_s}}(n_{i_s}-1)/2)}\cdot \left(\chi_{l_{i_s}}(-1) g(\chi_{l_{i_s}},\psi)\cdot
% \det(\Frob_{t_{l_{i_s}}}, F_{l_{i_s}})\right)^{n_{i_s}}\quad \textrm{if} \quad \chi_{l_{i_s}} \neq \1\\
% (-q)^{l\cdot n(n-1)/2}\det(-\Frob_s,F_{i_s}))^{n_{i_s}}\quad \textrm{if} \quad \chi_{l_{i_s}} =\1
%  \end{array}\right $$
% 
% \epsilon_0( T_{l_s},
% \cJ_{n_{i_s}}\otimes \Ind_{G_{s,l_{i_s}}}^{G_{l_s}}( V_{\chi_{i_s}}\otimes {F_{i_s}'}), d\pi)$$
%$$\begin{array{c}
%q^{l_{i_s}}\cdot (n_{i_s}}(n_{i_s}-1)/2)}\cdot \left(\chi_{l_{i_s}}(-1) g(\chi_{l_{i_s}},\psi)\cdot
% \det(\Frob_{t_{l_{i_s}}}, F_{l_{i_s}})\right)^{n_{i_s}}\quad \textrm{if} \quad \chi_{l_{i_s}} \neq \1\\
% (-q)^{l\cdot n(n-1)/2}\det(-\Frob_s,F_{i_s}))^{n_{i_s}}\quad \textrm{if} \quad \chi_{l_{i_s}} =\1
%  \end{array}$$
  % determined by  Prop.~\ref{propindecomp}.
  \end{thm}
  
  \proof   Let $X=\PP^1_k.$ The  cohomology groups
 $H^i(X_{\bar{k}},j_*{F_y(\chi)})\,(i\neq 1)$ vanish by the irreducibility assumption on~$F.$  
Hence with the conventions of Section~\ref{secfrobetc} we obtain 
 $$\det(-\Frob_q,R\Gamma(X_{\overline{k}},j_*{F_y(\chi)}))^{-1}=
\det(-\Frob_q, H^1(X_{\overline{k}},j_*{F_y(\chi)}) ).$$ Laumon's product formula for 
local epsilon constants (Formula~\eqref{eql4}) implies therefore that 
$$ \det(-\Frob_q,H^1(X_{\overline{k}},j_*{F_y(\chi)}))=
q^{r(F)}\prod_{x\in |X|}\epsilon(X_{(x)},F|_{\eta_{x}, }, \omega_0|_{X_{(x)}}).$$ 

We next use arguments similar to \cite{Laumon}, Prop.~3.3.2:
Since the valuation of $\omega_0|_{X_{(x)}}$ is zero for $x\in |\AA^1|$ it follows from 
Prop.~\ref{reml1}(ii),(iv) that for all $x\in |\AA^1\setminus S|$ one has 
$$ \epsilon(X_{(x)},{F_y(\chi)}|_{\eta_{x}, }, \omega_0|_{X_{(x)}})=1$$ 
%By the same argument, 
%the  epsilon constant of  the inertia invariants ${F_y(\chi)}_{\overline{\eta}_s}^{I_s^t^t}$ is equal to $1$ for all $s\in S.$  
and therefore 

\begin{equation}\label{eqfirstprod1} \det(-\Frob_q,H^1(X_{\overline{k}},j_*{F_y(\chi)}))=
q^{r(F)}\epsilon(X_{(\infty)},{F_y(\chi)}_{{\eta}_\infty},\omega_0|_{X_{(\infty)}})\prod_{s\in S\cup y}\epsilon(X_{(s)},{F_y(\chi)}_{{\eta}_s},  \omega_0|_{X_{(s)}}).\end{equation}

The sheaf  ${F_y(\chi)}$ is smooth at $\infty$ since $F$ is in standard form w.r. to $\chi.$ 
Since the valuation of $\omega_0|_{X_{(\infty)}}$ is equal to $-2$ 
it  follows from Prop.~\ref{reml1}(ii),(iv),(v) that 
\begin{equation}\label{eqfirstprod2} \epsilon(X_{(\infty)}, {F_y(\chi)}_{{\eta}_\infty}, \omega_0|_{X_{(\infty)}})
=  q^{-2r(F)}\cdot \det(\Frob_\infty,{F_y(\chi)}_{\bar{\infty}})^{-2}.\end{equation}

  To determine the  the other local epsilon constants, let $k\{\pi\}$ the henselization of $k[\pi]_{(\pi)}$ and $T=\Spec(k\{\pi\})=(T,\eta,t)$ the corresponding trait with uniformizer $\pi$ and tame 
  fundamental group $G^t.$  Let $s\in S\cup y$ with $k(s)=k_{l_s}=\FF_{q^{l_s}}.$ 
  Let $T_{l_s}:=T\otimes_k k_{l_s},$ with uniformizer $\pi_{l_s}$ be as in Prop.~\ref{reml1}(iii). 
  For each $s\in S$ we have a canonical uniformizer 
  $\pi_s$ of the trait $X_{(s)}$ as follows:
  for 
  $$ s\otimes \bar{k}=\coprod_{\iota\in \Hom_k(k(s),\bar{k})}s_\iota\quad \textrm{ define }\quad  \pi_s:=\prod_\iota(x-s_\iota).$$ 
  Then one obtains a finite \'etale $k$-morphism $f_s:X_{(s)}\to T$ induced by $\pi_s\mapsto
  \pi.$ 
  In this case we then have for any $G_s^t$-module $V$
  \begin{equation}\label{eqfirstprod3} \epsilon_0(X_{(s)},V, \omega_0|_{X_{(s)}})
  =\epsilon_0(T,\pi_{s,*}(V), d\pi),\end{equation}
  since we have $\omega_0|_{X_{(s)}}=f_s^*(d\pi)$ because  $f_s$ is basically the localization of the  inclusion $X_{(s)}\to X$  
  (Eq.~\eqref{eqfirstprod3} itself can be deduced using the last remark in Prop.~\ref{remdeligneaxioms}(iv) 
  as in \cite{Laumon}~(3.5.2), (3.5.3.1)). 
  Moreover, by Prop.~\ref{reml1}(iii), 
  \begin{equation}\label{eqfirstprod33}\epsilon_0(T,\pi_{s,*}(V), d\pi) =\epsilon_0(T_{l_s},V|_{T_{l_s}}, d\pi_{l_s})\,.\end{equation}
  
  This implies  one the one hand that 
  \begin{equation}\label{eqfirstprod5}\epsilon_0(\PP^1_{(y)},
 {F_y(\chi)}_{{\eta}_s}/{F_y(\chi)}_{{\eta}_s}^{I_s^t}, \omega_0|_{\PP^1_{(y)}})=
 (\chi(-1)g(\chi,\psi))^{r(F)}\cdot \det(\Frob_y,{F_y(\chi)}_{\bar{\eta}_y})\end{equation}
  by 
 Prop.~\ref{propindecomp} (or directly from Prop.~\ref{remdeligneaxioms} and \ref{reml1}). 

 On the other hand, by the same argument as used for Eq.~\eqref{eqfirstprod1}, 
the  local epsilon constant of  the inertia invariants ${F_y(\chi)}_{\overline{\eta}_s}^{I_s^t}$ is equal to $1$ for all $s\in S.$  
  Using the multiplicativity of $\epsilon_0$ on short exact sequences
%   together 
%  and a  base change to $T_{l_{i_s}}$ for each $i_s=1,\ldots,r_s,$ 
  together with 
   the last statement in Prop. \ref{propindecomp} one sees that 
  \begin{equation}\label{eqfirstprod6}
  \epsilon(X_{(s)},{F_y(\chi)}_{{\eta}_s},  \omega_0|_{X_{(s)}})=\prod_{i_s=1}^{r_s}
 \epsilon_0( \PP^1_{({s})},
 \cJ_{n_{i_s}}\otimes \Ind_{G_{s,l_{i_s}}}^{G_s}( V_{\chi_{i_s}}\otimes {F_{i_s}'}), \omega_0|_{\PP^1_{({s})}}) \,.  \end{equation} 
 Combining \eqref{eqfirstprod3} with \eqref{eqfirstprod33} and using Prop.~\ref{propindecomp} we obtain for $i_s=1,\ldots,r_s$
  \begin{equation}\label{eqfirstprod7}
 \epsilon_0( X_{({s})},
 \cJ_{n_{i_s}}\otimes \Ind_{G_{s,l_{i_s}}}^{G_s}( V_{\chi_{i_s}}\otimes {F_{i_s}'}), \omega_0|_{\PP^1_{({s})}}) = \end{equation}
 $$\left\{
 \begin{array}{ll}
q^{l_{i_s}\cdot n_{i_s}(n_{i_s}-1)/2)}
\cdot (\chi_{l_{i_s}}(-1) g(\chi_{l_{i_s}},\psi)\cdot
 \det(\Frob_{t_{l_{i_s}}}, F_{l_{i_s}}))^{n_{i_s}}& \quad \textrm{if}\quad  \chi_{l_{i_s}} \neq \1 \\
 (-q)^{l\cdot n_{i_s}(n_{i_s}-1)/2}\cdot \det(-\Frob_s,F_{i_s}))^{n_{i_s}}& \quad \textrm{if}\quad  \chi_{l_{i_s}} =\1\,,
  \end{array}\right.\,.   $$
 By combining the formulas \eqref{eqfirstprod1}, \eqref{eqfirstprod2}, \eqref{eqfirstprod6},
 and \eqref{eqfirstprod7}, one obtains the claimed result. 
%
%  Hence
%$$ \det(-\Frob_q,R\Gamma(X_{\overline{k}},j_*{F_y(\chi)}))^{-1}=
%q^{-r(F)}\cdot \det(\Frob_\infty,{F_y(\chi)}_{\bar{\infty}})^{-2}\cdot \prod_{s\in S\cup y}\epsilon_0(X_{(s)},{F_y(\chi)}_{{\eta}_s}/{F_y(\chi)}_{{\eta}_s}^{I_s}, \omega_0|_{X_{(s)}}),$$ implying the first equation.  
%
% The second equation is implied by Prop.~\ref{reml1}(ii),(iii), and (iv).
%  The third  equation follows again from Prop.~\ref{reml1}. 
%  
 % To see the last equality, we argue as follows (see \cite{Laumon}~(3.4.1.1)): 

%  Moreover,  for any $G_s^t$-module $V,$ 
%  $$ \epsilon_0(T,\Ind_{G_s^t}^{G^t}(V), -d\pi)=\epsilon_0(T,\Ind_{G_s^t}^{G^t}(V), -d\pi)
%  Therefore, 
%  $$  \epsilon_0(X_{(s)},\Ind_{G_{s,l_{i_s}}^t}^{G_s^t}( V_{\chi_{i_s}}\otimes F_{i_s}'), \omega_0|_{X_{(s)}})
%  Moreover, 
%  
%  see \cite{Laumon}, (3.5.3.1)

 %  using Prop.~\ref{propindecomp}.
   \Endproof

\begin{cor}\label{cormcchidet} Under the above assumptions,
\begin{eqnarray}
  \det(-\Frob_q,\MC_\chi(j_*(F)[1])_{\bar{y}})&=&\det(-\Frob_q,H^1(\PP^1_{\overline{k}},j_*{F_y(\chi)}))^{-1} \nonumber \\
  &=& q^{-r(F)}
 \det(\Frob_{\infty}, {F_y(\chi)}_{\bar{\infty}})^{-2} 
\prod_{s\in S\cup y}\epsilon(\PP^1_{(s)},
 F_{y,\overline{\eta}_{s}}, \omega_0|_{\PP^1_{(s)}}),\nonumber \end{eqnarray}
as in Thm.~\ref{thml2}. 
\end{cor}

\proof  Only the first equality has to shown. Using the definition
 of $\MC_\chi(F)$ (see  Eq.~\eqref{eqmcchidef} and Lemma~\ref{lemmiddle}) and 
 Prop.~\ref{rembasicprops}(iii)
 one sees that 
 $$\MC_\chi(j_*(F)[1])_{\bar{y}}=H^1(\PP^1_{\overline{k}},j_*{F_y(\chi)})[1] .$$  
 \Endproof

\begin{thm}\label{thml1} Let $k=\FF_q,$ be a finite field of odd order,
 let $U\stackrel{j}{\hookrightarrow} \AA^1_k$  be a dense  open subscheme, and let 
 $S=\AA^1\setminus U.$ Let $V=j_*F[1]\in {\rm Fourier}(\AA^1,\bQl)$  be an irreducible nonconstant 
tame middle extension sheaf, where $F$ is smooth on $U$ as above. 
%and in standard form with respect to the quadratic character $-\1:k^\times \to \bQl.$  
Assume that  $F$ satisfies the  the following conditions: 
\begin{enumerate}
 %\item {\rm (At most quadratic geometric determinant) }
 %The geometric determinant of the monodromy  of $V$ is either trivial
 %of given by a quadratic character of $\pi_1(U_{\overline{k}}).$ 
% \item The local geometric monodromy of $F$ at $\infty$ is scalar, given by the quadratic character $-\1:k^\times \to \overline{\mathbb{\QQ}}_\ell^\times,$
% but $F$ is not geometrically isomorphic to $ V_{-\1}.$  
 \item The sheaf $F$ is in standard form with respect to the quadratic character $-\1:k^\times \to \bQl.$  
  \item  % {\rm (Geometric rationality)} 
   The $I_s^t$-module $F_{\bar{\eta}_s}$  
   is self-dual for all $s\in S.$ 
   \item %{\rm (At most quadratic determinant)}
For each  $y\in |\AA^1_k|$  there exists an integer $m$ so that
     $$\det(\Frob_y,(j_*F)_{\overline{y}})=\pm  q^m.$$
%     \item For each $s\in S$ the  $1$-eigenspace of $ {F_y(\chi)}_{\overline{\eta}_s}/{F_y(\chi)}_{\overline{\eta}_s}^{I_s^t}$ is either trivial or it  decomposes into
%     (not necessarily indecomposable) summands of the form $\cJ_{k'}\otimes F'\,(k'\in \NN)$
%     such that $F'$ is unramified at $s$ and $$\det(\Frob_s,F')=\pm q^{n'}\quad \textrm{for some}\quad  n'\in \NN.$$ 
 \end{enumerate} 
Then $\MC_{-\1}(V)=j_*G[1]\in {\rm Fourier}(\AA^1,\bQl)$ with $G$ smooth on $U$ and 
$G$  satisfies the conditions (i)--(iii) with $F$ replaced by $G.$
%\begin{enumerate}
% %\item {\rm (At most quadratic geometric determinant) }
% %The geometric determinant of the monodromy  of $V$ is either trivial
% %of given by a quadratic character of $\pi_1(U_{\overline{k}}).$ 
%% \item The local geometric monodromy of $F$ at $\infty$ is scalar, given by the quadratic character $-\1:k^\times \to \overline{\mathbb{\QQ}}_\ell^\times,$
%% but $F$ is not geometrically isomorphic to $ V_{-\1}.$  
%  \item % {\rm (Geometric rationality)} 
%   The $I_s^t$-module $G_{\bar{\eta}_s}$ 
%   is self-dual for all $s\in S.$ 
%   \item %{\rm (At most quadratic determinant)}
%For each  $y\in |\AA^1\setminus S|$  there exists an integer $m'$ so that
%     $$\det(\Frob_y,(j_*G)_{\overline{y}})=\pm  q^{m'}.$$
%     \item[(iii')] For each $s\in S$ the  $-1$-eigenspace of $ G_{\overline{\eta}_s}/G_{\overline{\eta}_s}^{I_s}$ is either trivial or 
%     it  decomposes into
%     (not necessarily indecomposable) summands of the form $\cJ_{k''}\otimes G'\,(k''\in \NN)$
%     such that $G'$ is unramified at $s$ and $$\det(\Frob_s,G')=\pm q^{n''}\quad \textrm{for some}\quad  n''\in \NN.$$ 
% \end{enumerate} 
% Moreover, if $F$ satisfies the respective condition (iii'), then $G$ satisfies the 
% respective condition (iii). 
  \end{thm}
   
 \proof  It follows from Rem.~\ref{remstandardform} that $G$ is again in standard form 
 w.r. to $-\1$
and hence  fulfils the respective condition (i).
  The effect of $\MC_{-\1}$ on the geometric local monodromy (\cite{Katz96},~Thm.~3.3.5 and Cor.~3.3.6; cf.~Thm.~\ref{thmThom})
 implies that also (ii) holds for $G.$ 
 
% Moreover,   Thm.~\ref{thmThom} implies that the condition
% (iii) for $F$ implies condition (iii') for $G$ (and also the last statement of the theorem). 
%  
 The self-duality assumption in (ii) for $F$  implies that the contribution of the 
  Gau{\ss}-sums to the determinant of $G$ (by local epsilon constants as in  Thm.~\ref{thml2} and Cor.~\ref{cormcchidet}) which arise from nontrivial characters ${\chi}_{i_s}\neq -\1$
 cancel out to $\pm q^{k_1}\, (k_1\in \ZZ)$ because also ${\chi}_{i_s}^{-1}$ appears in
 this local monodromy (using Rem.~\ref{remgaussformula}). 
 
 The product formula for the (determinant of the) geometric monodromy implies that 
 the character $-\1$ appears in an even number. Hence also the contribution of these 
  Gau{\ss}-sums amounts to $\pm q^{k_2}\, (k_2\in \ZZ),$ again by Rem.~\ref{remgaussformula}. 
 
 It follows from condition (iii) for $F$  that for each $s\in S,$
 $$ \det(\Frob_s, F_y(-\1)_{{\eta}_s}/F_y(-\1)_{{\eta}_s}^{I_s^t})=\pm q^{k_3}$$
 for some $k_3\in \ZZ$ since tensoring by the quadratic Kummer sheaf 
 $\cL_{-\1}$ does not change this property.  
 Hence, by the usual rules for the behaviour of determinants under induction, 
  $$\prod_{i_s=1}^{r_s}\det(-\Frob_s,F_{i_s}')^{n_{i_s}}=\pm q^{k_4}$$ 
   for some $k_4\in \ZZ,$ where $F_{i_s}'$ is as in Thm.~\ref{thml2}.  
% 
% the trivial $I_s^t^t$-characters in the local monodromy (arising from purely unipotent indecomposable summands in $ F_{\overline{\eta}_s}/F_{\overline{\eta}_s}^{I_s^t}$)   amounts 
% to $\pm q^{k_3}\, (k_3\in \ZZ)$ by assumption~(iii) and Prop.~\ref{propindecomp}.
 Finally, the contribution of 
$ \det(\Frob_{\infty}, {F_y(\chi)}_{\bar{\infty}})^{-2} $ to the determinant of $G$ is $\pm q^{k_5}\, (k_5\in \ZZ)$ by assumption (iii) for $F$ and Cebotarev's density theorem. Via Thm.~\ref{thml2} this proves that the condition (iii) holds for $G$ for any closed point 
$y\in |\AA^1\setminus  S|.$

One has  $ \det(\Frob_s, {F}_{{\eta}_s})=\pm q^{k_6}$
for some $k_6\in \ZZ$ by Cebotarev's density theorem ($s\in S$).
It follows therefore from condition (iii) for $F$  together with ${F}_{{\eta}_s}^{I_s^t}\simeq 
(j_*F)_{\bar{s}}$ that 
 $$ \det(\Frob_s, {F}_{{\eta}_s}/{F}_{{\eta}_s}^{I_s^t})=\pm q^{k_7}\quad \textrm{for some}
 \quad k_7\in \ZZ. $$  Hence Thm.~\ref{thmThom} implies that 
$$ \det(\Frob_s, {G}_{{\eta}_s}/{G}_{{\eta}_s}^{I_s^t})=\pm q^{k_8}\quad \textrm{for some}
 \quad k_8\in \ZZ, $$ implying that the condition (iii) holds for $G$ also at each point $s\in S$
 (again by Cebotarev's theorem). 
 \Endproof

 \bibliographystyle{plain}
                    \bibliography{p1}

\end{document}